\documentclass{nm1}
\usepackage{bm}
\usepackage{charter}
\usepackage[charter]{mathdesign}
\usepackage{ulem}
\usepackage{epstopdf}
\usepackage{graphicx}
\usepackage{caption}
\DeclareCaptionLabelFormat{mylabel}{#1 #2.\hspace{1.5ex}}
\usepackage{subcaption}
\usepackage{enumerate}
\usepackage{booktabs}
\usepackage[percent]{overpic}

\usepackage{url}
\usepackage{algorithm}
\usepackage[noend]{algpseudocode}
\def\R{\mathbb{R}}

\usepackage{color}

\newcommand{\blue}[1]{{\color{blue} #1}}

\begin{document}
	
	
	\markboth{S. Li, L. Ling, \& et al.}{Enhancing RBF-FD Efficiency for Highly Non-Uniform Node Distributions via Adaptivity}
	\title{Enhancing RBF-FD Efficiency for Highly Non-Uniform Node Distributions via Adaptivity}
	

	%
	%
	%
	
	\author[AUTHOR1 and AUTHOR2]{Siqing Li\affil{1}, Leevan Ling\affil{2} \comma \corrauth, Xin Liu\affil{3} , Pankaj K. Mishra\affil{4}, Mrinal K. Sen\affil{5}, and Jing Zhang\affil{1}}
	\address{\affilnum{1}\ College of mathematics, Taiyuan University of Technology, Taiyuan,  China. \\
		\affilnum{2}\ Department of Mathematics, Hong Kong Baptist University, Hong Kong.
		\\
		\affilnum{3}\  \blue{Institute of Advanced Technology, University of Science and Technology of China, Hefei, China.}
		\\
		\affilnum{4}\  Geological Survey of Finland.
		\\
		\affilnum{5}\ \blue{ University of Texas at Austin, USA.}
	}
	
	\emails{{\tt lisiqing@tyut.edu.cn} (S. Li)
	}

	\begin{abstract}
		Radial basis function generated finite-difference (RBF-FD) methods have recently gained popularity due to their flexibility with irregular node distributions. However, the convergence theories in the literature, when applied to nonuniform node distributions, require shrinking fill distance and do not take advantage of areas with high data density. Non-adaptive approach using same stencil size and degree of appended polynomial will have higher local accuracy at high density region, but has no effect on the overall order of convergence and could be a waste of computational power.
		This work proposes an adaptive RBF-FD method that utilizes the local data density to achieve a desirable order accuracy. By performing polynomial refinement and using adaptive stencil size based on data density, the adaptive RBF-FD method yields differentiation matrices with higher sparsity while achieving the same user-specified convergence order for nonuniform point distributions. This allows the method to better leverage regions with higher node density, maintaining both accuracy and efficiency compared to standard non-adaptive RBF-FD methods.
	\end{abstract}
	
	
	\maketitle
	
	\section{Introduction}
	\label{sec1}
	In the past two decades, there has been important progress in developing adaptive mesh methods for PDEs. Mesh adaptivity is usually of two types in form: local mesh refinement and moving mesh method.
	
	Radial basis functions (RBFs) have been a popular choice in the development of kernel-based meshless methods for solving partial differential equations (PDEs) numerically. Besides collocation methods, the localized RBF-FD method has gained popularity in recent years due to its many advantages, including numerical stability on irregular node layouts, high computational speed and accuracy, easy local adaptive refinement, and excellent opportunities for large-scale parallel computing.
	
	The idea of RBF and kernel-based differentiation can be traced back to \cite{Tolstykh2003}, and it was formally introduced as the RBF-FD in \cite{wright2003}. Since then, a significant amount of research has been dedicated towards the robust development of the RBF-FD method \cite{wright2006,bayona2010,fornberg2011,fornberg2015,Flyer2016role,Bayona2017role,shankar2017,santos2017,petras2018rbf,Mishra2018}, as well as its application to various problems in science and engineering \cite{chandhini2007local,chinchapatnam2009compact,flyer2012guide,shankar2015radial,martin2015seismic,martin2017using,mishra2017frequency,kindelan2018frequency,slak2018refined,martin2018improved}. In addition to the RBF-FD method, other collocation methods based on radial basis functions have also been proposed. For example, a global radial basis function collocation method in \cite{arler2005ARB} was successfully applied to solve a computational fluid dynamic problem, and local RBF collocation methods were used to solve the diffusion problem in \cite{SARLER20061269} and Hamiltonian PDEs in \cite{zhang2021}.

	The RBF-FD method is advantageous since it work with scattered nodes, allowing for stencils with different configurations and overcoming the fixed grid/element limitation of conventional numerical methods. Unlike global meshless methods, RBF-FD computes weights locally using RBFs expanded at a fixed number of nearest nodes. Once weights at each node are computed, they can be stored and used for next-step computation, making weight computation a pre-processing step in solving time-dependent PDEs. Furthermore, weight computations at each node are independent processes, making RBF-FD a desirable method for parallel computing.
	
	It has been demonstrated \cite{Flyer2016role,Bayona2017role} that combining  polynomial basis with polyharmonic spline radial basis functions (PHS+Poly) in the RBF-FD formulation leads to considerable improvements in robustness. Key benefits of the PHS+Poly approach include:
	\begin{enumerate}
		\item[1.]  It is free of shape parameter, simplifying the formulation and eliminating the need for fine-tuning.
		\item[2.] It is efficient compared to stable RBF-FD formulations based on infinitely smooth RBFs \cite{santos2017}.
		\item[3.] This method ensures accuracy near boundaries without ghost-nodes where stencils become highly one-sided \cite{BAYONA2019}.
		\item[4.] It has the potential to maintain accuracy for large and sparse linear systems.
		\item[5.] The convergence order depends mainly on the augmented polynomial degree, which determines the stencil size.
	\end{enumerate}
	
	Existing convergence theories for scattered nodes require shrinking fill distance and fail to leverage regions of high node density. Numerically, non-adaptive methods using uniform stencil size and polynomial degree may exhibit loss of accuracy in low-density regions while adding unnecessary complexity in high-density regions, decreasing computational efficiency.
	Rather than a fixed stencil size and augmented polynomial degree, we propose an adaptive RBF-FD algorithm that allows the user to define a global convergence order  with respect to the total number of nodes in the domain. The method can efficiently obtain the user-defined convergence order in a highly nonuniform node-layout with a significantly large mesh ratio. Note that the mesh ratio is defined as the fill distance over minimum separating distance, which plays a role similar to mesh ratio in finite element methods.
	
	The rest of the paper is organized as follows. Section~\ref{sec:rbffd} discusses the general formulation of RBF-FD and polynomial augmentation, along with some insights for its stable implementation.  In Section~\ref{sec:adaprbffd}, we provide a precise definition of  the global convergence order  in a nonuniform node-layout.  We also present a method that connects the global convergence order to the required degree of polynomial at a local scale, depending on the local fill distance. This is followed by a discussion of the weight computation approach through the adaptive RBF-FD method. In Section~\ref{sec:num}, numerical examples are illustrated to show the advantages of the proposed method. Finally, section~\ref{sec:con} is the conclusions drawn from our study.

	\section{RBF-FD approximation} \label{sec:rbffd}
	\begin{figure}[t]
		\centering
		\includegraphics[scale=0.6]{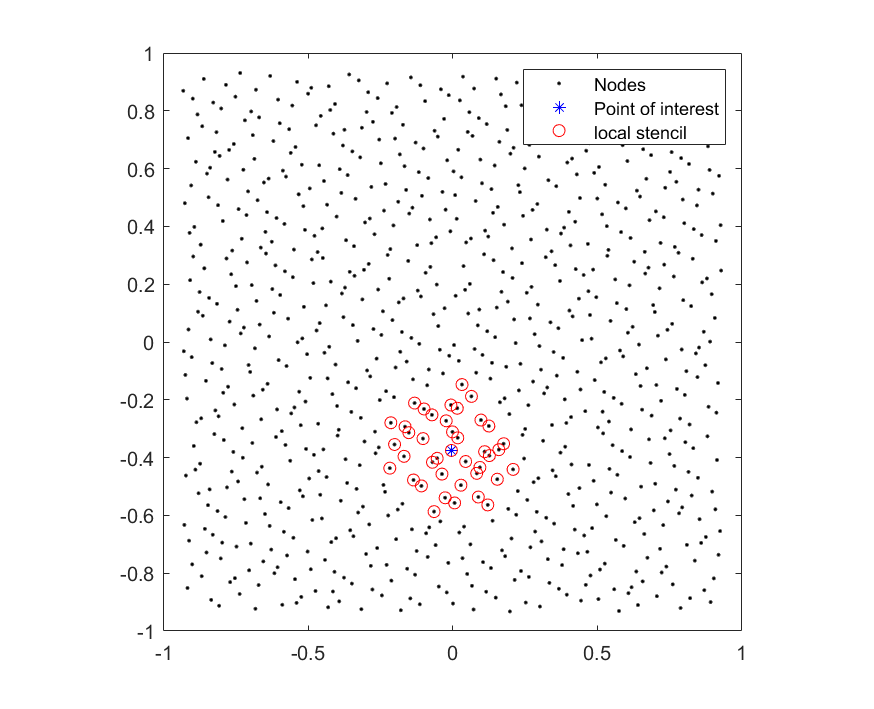}
		\caption{   A prototype figure showing centers and local stencil.
		}
		\label{fig:scatterdata}
	\end{figure}
	
	We first review the basic RBF-FD method.
	Let $\Omega \subset \mathbb{R}^d$ be a computational domain with a set of distinct scattered nodes $X:=\{\bm{x}_1, \bm{x}_2, ..., \bm{x}_N\} \subset \Omega$.
	Given a linear differential operator $\mathcal{L}$, we aim to approximate the values of $\mathcal{L}u$ at a center/point of interest $\bm{x}^c \in X$ based on the nodal values of $u$ at some stencil or set of neighboring nodes of $\bm{x}^c$  in $X$.  Figure \ref{fig:scatterdata} shows the prototype of the scatter data distribution and local stencil for the point of interest.

	Let the stencil of $\bm{x}_c$ be denoted by $X_c:=\{\bm{x}_1^c,\ldots, \bm{x}_{n_s}^c\}\subseteq X$, where $n_s^c$ is the stencil size of $\bm{x}_c$.
	Then, the value of $\mathcal{L}u(\bm{x}_c)\in\R$ can be approximated by  as:
	\begin{equation}
		\mathcal{L}u(\bm{x}_c) \approx \sum_{k=1}^{n_s^c} w_k^c u\left(\bm{x}_k^c \right)=\bm{w}_c u|_{X_c},
		\label{eq:interp}
	\end{equation}
	where $u|_{X_c}:= [ u(\bm{x}_1^c),\ldots, u(\bm{x}_{n_s}^c) ]^T $ is the nodal values of $u$ evaluated at stencil $X_c$.
	We need to compute the weights $\bm{w}_c:=[w_1^c, w_2^c, \ldots, w_N^c]^T$ at $\bm{x}_c$ for this finite difference approximation.

	In RBF-FD method, interpolants were used as surrogate. Using some symmetric positive definite (SPD) radial basis function  $\phi : [0,\infty) \rightarrow \mathbb{R}$, the RBF surrogate interpolates the data at $(X_c,u|_{X_c})\in \R^d\times \R$ while keeping function values $u|_{X_c}$ as unknown.
	The  RBF-FD weights $\bm{w}_c$ are computed by applying $\mathcal{L}$ to the surrogate and, then, evaluating  at $\bm{x}_c$. This is equivalent to solving:
	\begin{equation}
		\underbrace{\left[ \phi( \|\bm{x}_i^c-\bm{x}_j^c\|) \right]_{1\leq i,j \leq n_s^c} }_{\displaystyle =:A_c}\bm{w}_c
		= \underbrace{\left[ \mathcal{L} \phi( \|\bm{x}-\bm{x}_i^c\|)_{| \bm{x}=\bm{x}_c} \right]_{1\leq i \leq n_s^c }}_{\displaystyle =: [\mathcal{L}\phi]_{|\bm{x}_c}}.
		\label{eq:linearsys}
	\end{equation}
	Since $\phi$ is SPD, $A_c$ is SPD and (\ref{eq:linearsys}) is uniquely solvable.
	It is shown in \cite{Schaback-comput.oolcompline:14} that (\ref{eq:linearsys}) is error optimal in the native space norm corresponding to the radial basis function $\phi$.  We refer readers to the article for details.
	
	\subsection{RBF-FD using PHS+Poly}
	In contrast to global meshless methods, the local nature of RBF-FD rules out exponential convergence for approximating functions in certain native spaces. Even with infinitely smooth RBFs, the convergence is limited for RBF-FD by stencil sizes. Motivated by the discussion in \cite{Schaback-ErroAnalNodaMesh:17}, we use a conditionally positive definite kernel with finite order of smoothness in the RBF-FD method. For example, we and use a conditionally SPD polyharmonic spline (PHS) kernel:
	$$\phi (r) = r^{2m+1}, \  m \in N^+.$$
	To ensure unique solvability, we  augment polynomial basis of sufficient degrees  to ensure exact reproduction of low order polynomials. To compute RBF-FD weights in  (\ref{eq:interp}) by PHS+Poly,   the surrogate $I(\mathbf{x})$ now uses RBFs and polynomial basis up to degree $p$ as basis functions:
	\begin{equation}
		I(\mathbf{x}) \approx \sum_{k=1}^{n_s^c} \lambda_k \phi\left( \parallel\bm{x} -\bm{x}_k^c\parallel\right) + \sum_{{l}=1}^{n_p^c}\gamma_{{l}} p_{l}(\bm{x}),\quad \bm{x}\in \Omega,
	\end{equation}
	subject to interpolation conditions at data $(X_c,u|_{X_c})\in \R^d\times \R$ and  constraints
	\begin{equation}
		\sum_{k=1}^{n_s^c} \gamma_l p_l(\bm{x}_k^c) =0, \quad l=1,2,3,...,n_p^c,
	\end{equation}
	where  $n_p^c$ is the number of augmented multivariate polynomial basis and $n_s^c$ is again the stencil-size of $\bm{x}_c$  .

	The counterpart to the linear system (\ref{eq:linearsys}) for RBF-FD using PHS+Poly is given by
	\begin{equation}\label{eq:weights}
		\begin{bmatrix}
			A_c & P_c \\
			P_c^T & O
		\end{bmatrix}
		\begin{bmatrix}
			\bm{w}\\
			\bm{w}_e
		\end{bmatrix} =
		\begin{bmatrix}
			[\mathcal{L}\phi]_{|\bm{x}_c} \\
			[\mathcal{L}{p}]_{|\bm{x}_c}
		\end{bmatrix},
	\end{equation}
which can be interpreted as an equality-constrained quadratic programming problem (see \cite{Bayona2017role} for more details on this).
	In \eqref{eq:weights}, matrices $A_c$ and $[\mathcal{L}\phi]$ are defined as in (\ref{eq:linearsys}).  The $n_s^c \times n_p^c$ matrix $P_c$ and $n_p^c \times 1$ vector $[\mathcal{L}{p}]_{|\bm{x}_c}$ are given as
	\[
	P_c = [p_j(\bm{x}_i^c)]_{1\leq i \leq n_s^c,\, 1\leq j \leq n_p^c},
	\mbox{\quad and \quad}
	[\mathcal{L}{p}]_{|\bm{x}_c} = [ \mathcal{L}p_l(\bm{x})_{|\bm{x}=\bm{x}_c} ]_{1\leq l\leq n_p^c}.
	\]
	From here and on, we drop the super/subscript $c$ from  notations for simplicity, and all computations need to be done center by center for all $\bm{x}_c\in X$.
	The number of polynomial  basis $n_p$ up to degree $p$ in $d$-dimensional is given by
	\begin{equation}
		n_p = C(p+d,d) = \frac{(p+d)!}{p!d!}.
	\end{equation}
	Table \ref{tab:polterms} lists out polynomial terms for different degrees in some cases in 2D, see \cite{Flyer2016role} for more. For example, $n_p =3$ in 2D, corresponds to appending polynomial up to  first order.

	\begin{table}
		\centering
		\caption{The polynomial basis of different degree for 2D case }
		\begin{tabular}{ccc}
			\toprule
			Polynomial degree ($p$) & Polynomial basis & $n_p = C(p+2,2)=\frac{(p+2)!}{p!2!}$ \\
			\midrule
			0 & $1$ & 1 \\
			1 & $1,x,y$ & 3\\
			2 & $1,x,y,x^2,y^2,xy$ & 6\\
			3 & $1,x,y,x^2,y^2,xy,x^3,y^3,x^y,xy^2$ & 10 \\
			\bottomrule
		\end{tabular}
		\label{tab:polterms}
	\end{table}
	
\section{Adaptive PHS+Poly RBF-FD   methods }\label{sec:adaprbffd}
Adaptive methods are quite popular in the context of finite element methods (FEM), and  many  algorithms such as $h$-FEM, $p$-FEM, or $hp$-FEM were proposed. In order to meet a certain order of accuracy throughout the domain, the idea behind adaptivity is to vary parameters affecting the accuracy of the algorithm locally.
In this section, we aim to develop an adaptive PHS+Poly RBF-FD scheme that:
\begin{enumerate}
	\item[1.] Achieve the expected order of accuracy  in the domain, including evaluation points near boundaries where stencils are highly one-sided.
	\item[2.] Maintain accuracy and obtain linear systems with higher sparsity for dealing with large-scale problems.
\end{enumerate}

\subsection{Polynomial refinement and adaptive stencil size}

\begin{figure}[t]
	\centering
	\includegraphics[scale=0.75]{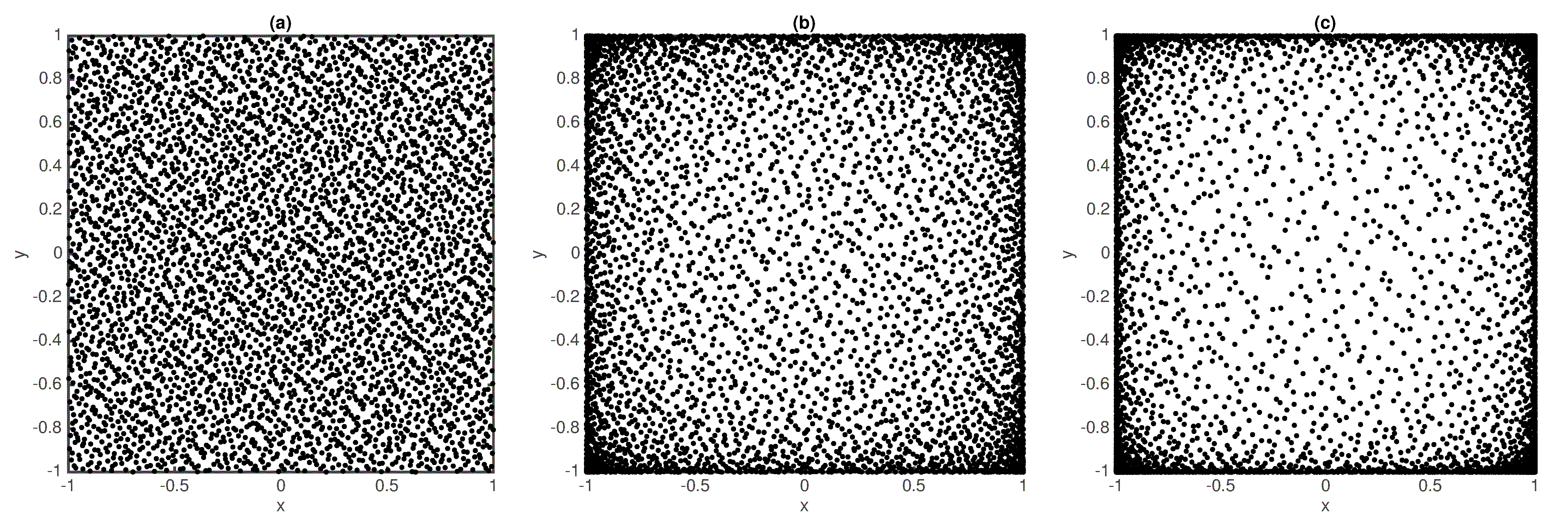}
	\caption{Node distribution under different mesh ratio with $N=4096$  quasi-uniform data points in unit square $[-1,1]^2$ with various mesh ratios: (a)  $\rho_X \approx  1$, (b) $\rho_X\approx 10^{2}$, and (c) $\rho_X \approx 10^{5}$.}
	\label{F:NT1nodes}
\end{figure}

Given a set of $N$ scattered data $X\subset\mathbb{R}^d$. Its fill distance and separation distance are given by
\begin{equation}\label{h}
	h_{X,\Omega} := h_X= \sup_{\bm{x}\in\Omega}\min_{\bm{x}_j\in X}\|\bm{x}-\bm{x}_j\|_{\ell_2(\R^d)},\quad q_X:= \frac12 \min_{i\neq j}\|\bm{x}_i-\bm{x}_j\|_{\ell_2(\R^d)},
\end{equation}
and $\rho_X = h_X/q_X \geq 1$ is the mesh ratio of $X$.
A quasi-uniform data set refers to a sequence of data sets where the mesh ratio $\rho_X$ remains bounded for all $N$ as $X$ refines with increasing $N$. The number of nodes $N$ typically scales as $N\sim q_X^{-d}$, where $\sim$ means proportional up to a constant. Figure \ref{F:NT1nodes} shows node distributions of $N=$ 4,096 under mesh ratios $\rho_X=$ 1, $10^2$, and $10^5$, respectively.

Data sets with large mesh ratios will be considered highly nonuniform, which is the target of our proposed adaptive approach. Nonuniform nodes are often associated with rapidly varying regions of a function or other physical quantities. 
For example, the adaptive node distributor in \texttt{NodeLab}\blue{ \cite{Mishra2019nodelab,FORNBERG2015531} }can generate a sequence of scattered data sets with large mesh ratios $\rho$ that are denser in user specified problematic spatial regions.
The folklore is that high point density will induce a smaller stencil radius even if the stencil size $n_s$ is fixed at a constant for all centers. Yet this could yield an inefficient algorithm, which we aim to improve upon.

In the RBF-FD method, for quasi-uniform node distribution, the convergence can be estimated in a straightforward manner. Suppose we want to approximate a $k$-th order differential operator $\mathcal{L}$ by PHS+Poly RBF-FD with appended polynomial up to degree $p$.  According to results in \cite{Flyer2016role,Bayona2017role}, for a quasi-uniform node distribution, the convergence order of a RBF-FD method is $\sim h_X^{p-k+1}$.
It has been found in \cite{Flyer2016role,Bayona2017role} that the convergence order of PHS+Poly  RBF-FD was governed by the degree of augmented polynomial $p$, while the smoothness order of the PHS  has a marginal effect on the accuracy and no influence on the order of convergence.

Let us first consider an ideal quasi-uniform point set $X$ of $N$ points in $\Omega$ (with minimal mesh ratio) that has an equispaced fill distance proportional to
\[
h_{e,\Omega} = \left( \frac{\text{Vol}(\Omega)}{N} \right)^{1/d}.
\]
Suppose we want to design an adaptive PHS+Poly RBF-FD  algorithm to approximate a $k$-th order differential operator $\mathcal{L}$ with user-specified $g^{th}$ order convergence. The approximation error is then expected to behave like $\varepsilon \sim h_{e,\Omega}^g$.
Moving to highly nonuniform nodes with large mesh ratio, we define global convergence with respect to $h_{e,\Omega}$ by formulating an adaptive RBF-FD approach with adaptive polynomial refinement.

For a general data set $X\subset \Omega$ with no constraint on its mesh ratio, let $X_c\subset X$ denote the stencil associated with the point of interest $\bm{x}^c \in X $,  and $\Omega_c = conv(X_c)$ be its associated local sub-domain.
Also let $\Omega_c:=conv(X_c)$ be its local sub-domain defined by convex hulls and $h_{X,\Omega_c}$ be the associated local fill distance.
If we employ a $g^{th}$-order RBF-FD scheme at this center, then the local error $\varepsilon_c$ will reduce as
\[
\varepsilon_c \sim h_{X}^g = (\rho_X q_X)^g \sim \rho_X^g N^{-g/d}.
\]
For uniform node with $\rho_X=1$, the last factor $N^{-g/d}\sim h_e^{g}$ is exactly what we expect from a $g^{th}$-order finite difference scheme.
In cases when $\rho_X\gg1$,  the asymptotic convergence rate remains order-$g$ but with a large leading constant $\rho_X^g$ that reduces the final accuracy in the error bound. To sum up, the large mesh ratio of highly nonuniform nodes weakens the approximation accuracy in terms of magnitude, but not order of convergence.

Since the (local) convergence order by PHS+Poly RBF-FD   algorithm is decided by the polynomial order and (local) fill distance $h_{X,\Omega_c}$,   we can obtain   $g^{th}$-convergence by adaptively varying the degree $p_c$ of augmented polynomial order that satisfies:
\begin{equation}\label{eq:cov1}
	h_{X,\Omega_c}^{p_c-k+1} = h_{e,\Omega}^g.
\end{equation}
In this context, non-adaptive standard approach takes $p_c-k+1=g$ for all centers $\bm{x}^c$.
Note that, by the definition of fill distance, at least one of these local fill distance $h_{X,\Omega_c}$ must coincide with the global fill distance $h_{X,\Omega}$, that acts as an  upper bound. This particular stencil is the limiting factor for the final accuracy.

Given the equispaced fill distance $h_{e,\Omega}$ (or $\text{Vol}(\Omega)$ and $N$), the local fill distance $h_{X,\Omega_i}$ for each point of interest, and  the order $k$ of the  differential operator $\mathcal{L}$,   the local polynomial orders for $g^{th}$ convergence can be determined from \eqref{eq:cov1} as:
\begin{equation}\label{eq:thumbrule}
	\mathbb{N}\ni p_{c} \geq g\frac{\log_{10}h_{e,\Omega}}{\log_{10}h_{X,\Omega_c}}+k-1.
\end{equation}
For quasi-random nodes with small mesh ratio,  we expect $h_{X,\Omega_i} \approx h_{e,\Omega}$ and $p_c$ in (\ref{eq:thumbrule}) simplifies to
\begin{equation}
	p = g+k-1 \qquad  \text{or} \qquad g = p-k+1.
	\label{eq:thumbrule2}
\end{equation}
Our adaptive scheme aligns with  existing non-adaptive methods \cite{Flyer2016role,Bayona2017role} in this case.  By adaptively determining $p_c$ from \eqref{eq:thumbrule}, we can achieve better accuracy for highly nonuniform nodes.

The stencil size $n_s$ also impacts the accuracy and convergence of RBF-FD. Studies have shown that when the polynomial order is sufficiently high, increasing $n_s$ beyond the number of polynomial basis functions $n_p$ does not significantly improve accuracy or convergence order, see  \cite{Flyer2016role,Bayona2017role}. For lower-order polynomials, increasing $n_s$ can improve accuracy but convergence order remains unchanged. This independence on stencil size is important for our adaptive scheme. In denser regions, we can use a smaller polynomial degree (and $n_s$) without significant losses. A stencil size of $n_s = n_p$ is suggested in \cite{Flyer2016role,Bayona2017role} for distributions with ghost nodes outside the boundary to avoid stagnation error. Using  PHS only without appended polynomials in RBF-FD, stagnation errors stem mainly from boundary errors. However, with polynomial augmentation and $n_s \geq 2n_p$, RBF-FD remains free from stagnation error.

\begin{figure}
	\centering
	\includegraphics[width=1\textwidth]{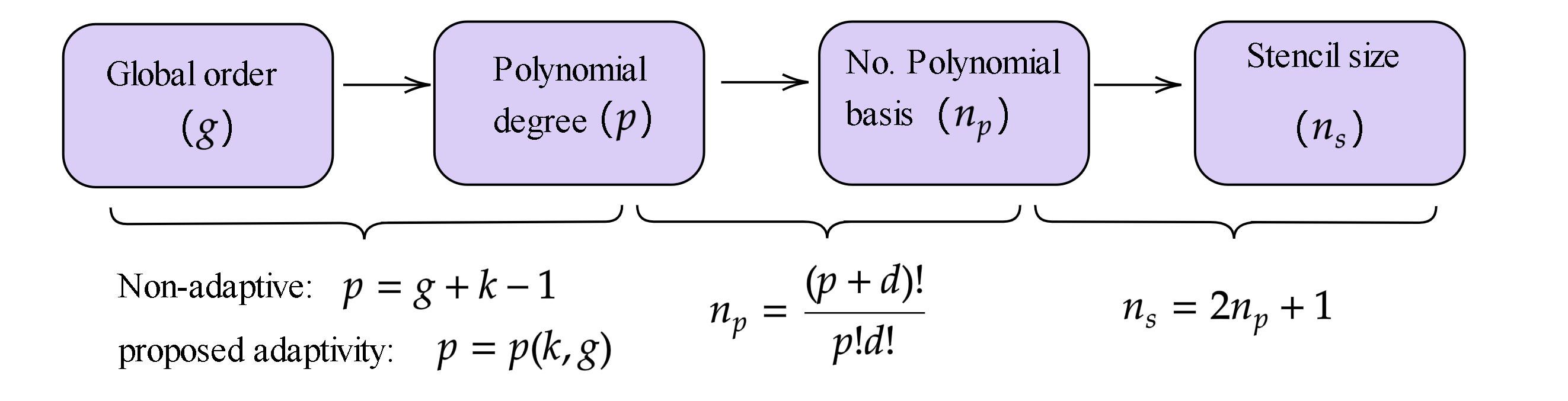}
	\caption{ The process of achieving a global order of convergence in an PHS+Poly RBF-FD  formulation. The order of convergence  $g$ input by the user.  The algorithm decides the degree of polynomial to be augmented, which also gives the number of polynomial terms required ($n_p$) and consequently the stencil size ($n_s$). }
	\label{F:g2p}
\end{figure}

Based on these findings, we propose:  after computing $p_c$ from (\ref{eq:thumbrule}) for a user-specified order of convergence $g$, we use $n_s = 2n_p+1$ as suggested in \cite{Flyer2016role,Bayona2017role} to ensure the solvability and the stability of our algorithm.
Since oversized stencil size  does not improve accuracy,   a smaller $n_s$  can be used in regions with smaller $h_{X, \Omega_i}$without significantly affecting accuracy or stability.
Figure \ref{F:g2p} illustrates our proposed adaptive PHS+poly RBF-FD scheme. Since the computational cost of the RBF-FD depends on the stencil-size $n_s$, using smaller local $n_s$ in denser regions can improve efficiency. The resulting FD-differentiation matrix also enjoys higher sparsity.
Algorithm 1 summarizes the RBF-FD weight computation using our adaptive PHS+Poly RBF-FD scheme.

\begin{algorithm}[t]
	\small
	\textsc{ {\textbf{Algorithm 1:} Adaptive PHS+Poly RBF-FD  method}}
	\begin{algorithmic}[1]
		\Procedure{Input}{$g\in \mathbb{N}$, $X\subset\Omega\subset\mathbb{R}^d$}
		\Comment{\scriptsize Input expected order and nodes}
		\State {Compute total number of nodes $N=|X|$}
		\State {Compute $h_{e,\Omega}= \left({\frac{Vol(\Omega)}{N}}\right)^{1/d}$}
		\State {Determine maximum stencil size} $n_{s,\max}$
		\For{\textsc{each point of interest $\bm{x}_c$} }
		\State {Find $n_{s,\max}$ -- nearest neighbors for each $\bm{x}_c$}
		\Comment{\scriptsize We have used knn-search}
		\State Find $h_{X,\Omega_c}$
		\Comment{\scriptsize  Estimate local fill-distance}
		\State $p_{i} = \lceil g\frac{\log_{10}h_{e,\Omega}}{\log_{10}h_{X,\Omega_c}}+k-1\rceil$
		\Comment{\scriptsize Polynomial degree required to maintain order $g$}
		\State $ n_p = C(p+d,d) = \frac{(p+d)!}{p!d!} $
		\Comment{\scriptsize Number of polynomial basis}
		\State $ n_s = 2n_p+1$
		\Comment{\scriptsize Stencil size just enough to support the polynomial $p$}
		\State $r = \sqrt{(x^2+y^2)}$
		\State $ A0 = r^m $
		\Comment{\scriptsize  The RBF interpolant}
		\State $ RHS0 = [rhs_L]$
		\Comment{\scriptsize  The RHS column in equation (\ref{eq:interp})}
		\State $ A = A0+Poly(p)$
		\Comment{\scriptsize   The RBF interpolant after augmenting polynomial of degree $p$}
		\State $ RHS = RHS0+Poly(p)$
		\Comment{\scriptsize Adding corresponding polynomial of degree $p$}
		\State $ w_L = A/RHS $
		\Comment{\scriptsize Solve the linear system for weights}
		\State $ w_L = w(1:n_s) $
		\Comment{\scriptsize Ignore the weights corresponding to polynomial augmentation}
		\EndFor
		\State \textbf{end}
		\EndProcedure
		\State \textbf{end}
	\end{algorithmic}
	\label{alg2}
\end{algorithm}

\section{Numerical experiments}\label{sec:num}
This section presents four numerical experiments demonstrating the performance of the proposed adaptive method. First, illustrative examples show the key concepts and behavior of the adaptive method. Then, a benchmark problem evaluates accuracy and convergence. Finally, real-world applications demonstrate effectiveness for both steady and time-dependent problems, showing the method's robustness, accuracy, and adaptivity.

\subsection{Example 1: Performance of the proposed adaptive method}
\begin{figure}
	\centering
	\centering
	\subcaptionbox{}{
		\includegraphics[width=6cm]{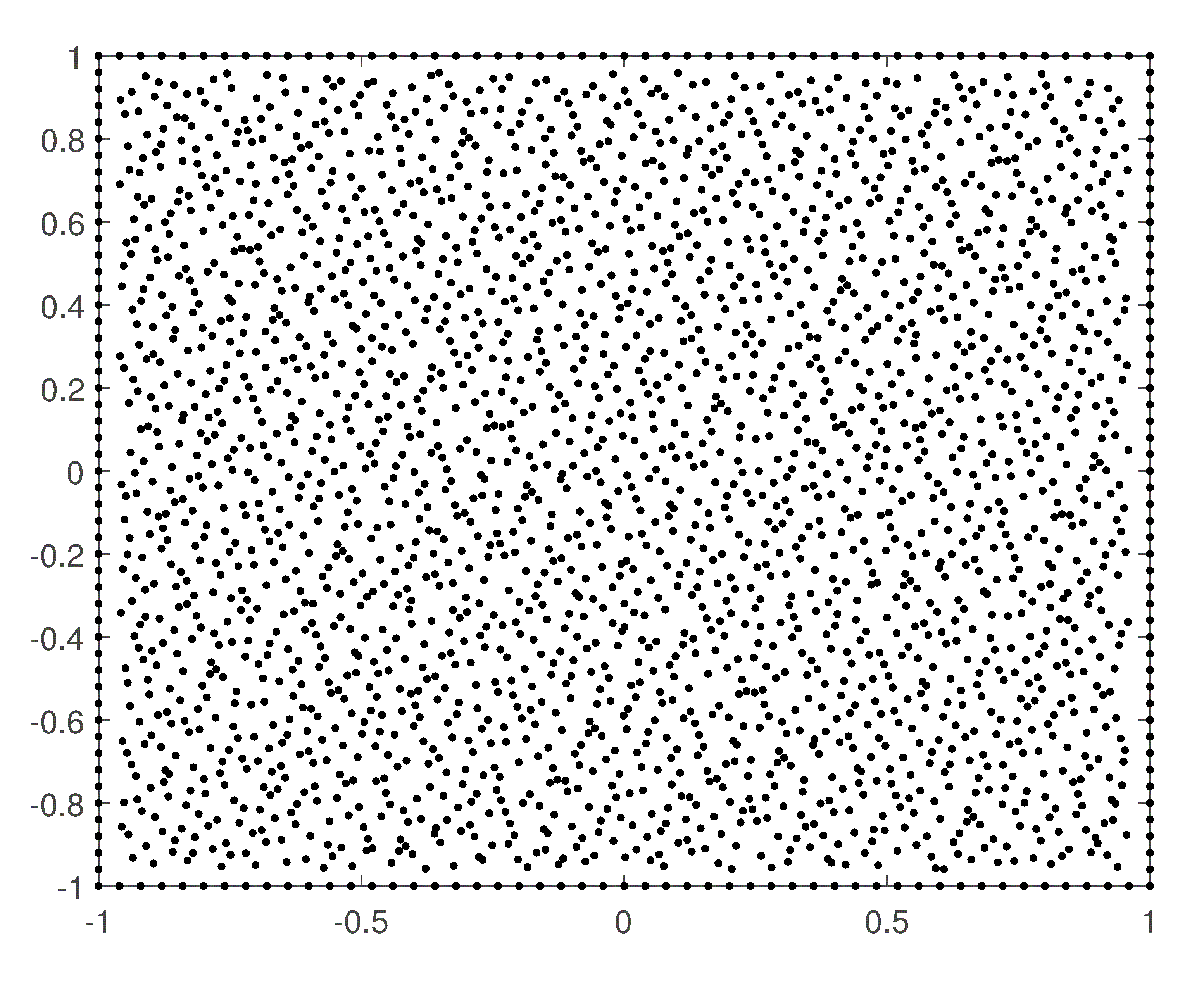}
	}
	\quad
	\subcaptionbox{}{
		\includegraphics[width=6cm]{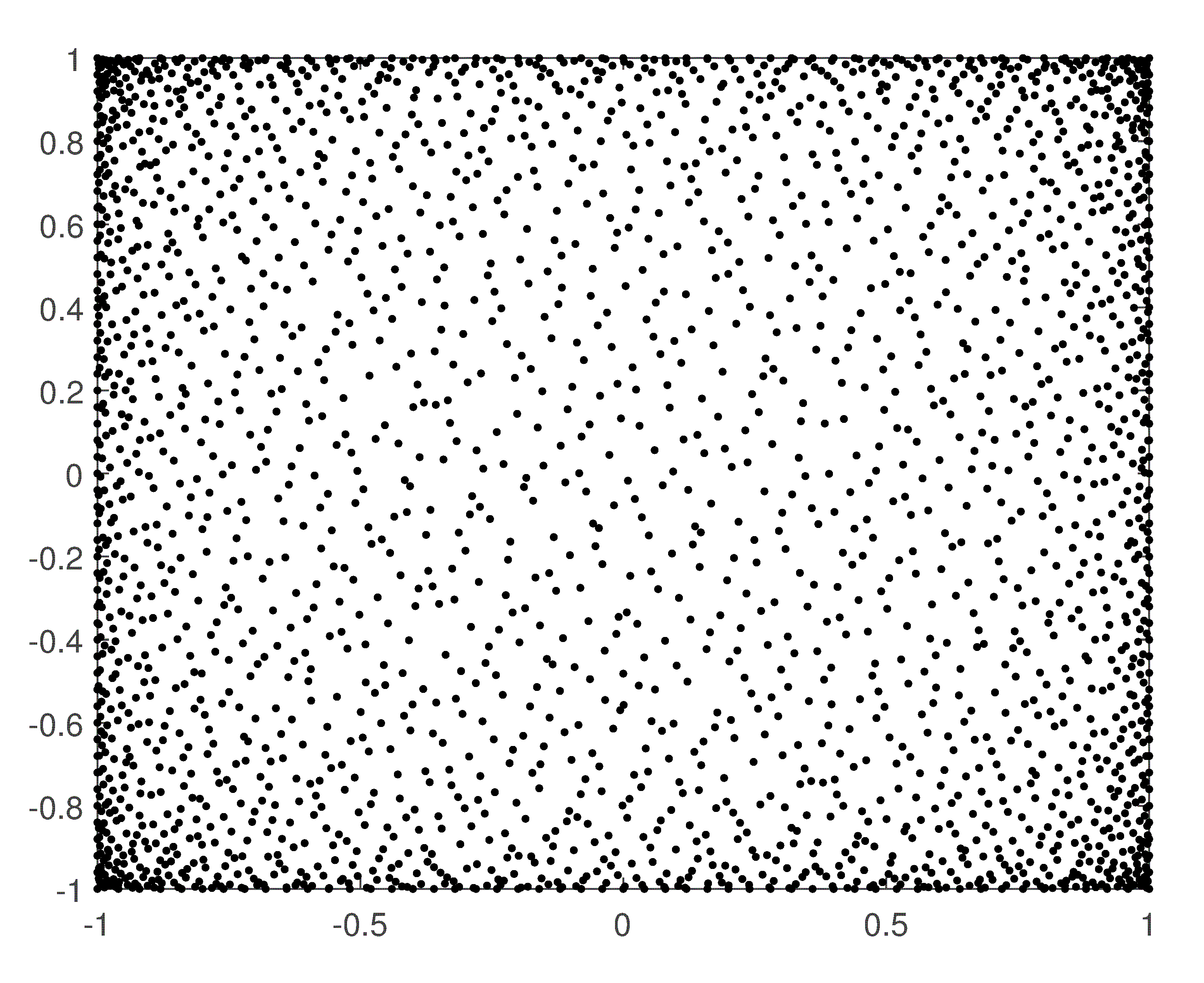}
	}
	\caption {Example 1:  { in $\Omega = [-1,1]^2$,  setting $N=2500$  (a) Halton points, and (b)} {nonuniform  nodes with large mesh ratio by sine-transform} .}
	\label{F:NT2Nodes}
\end{figure}

In the first test, we investigate the proposed adaptive PHS+Poly RBF-FD scheme by solving an elliptic partial differential equation in domain $\Omega = [-1,1]^2$ with  boundary  $\Gamma_1: \{(x,y)|-1\leq x \leq 1, y=1\}$ and $\Gamma_2:=\partial \Omega \setminus \Gamma_1$. The problem is set up as follows: find $u(x,y)$ such that
\begin{eqnarray*}
	\begin{aligned}
		-\nabla^2 u(x,y)  & = f(x, y),\qquad \qquad (x,y) \in \Omega,
	\end{aligned}
\end{eqnarray*}
subject to boundary conditions
\begin{eqnarray*}
	\begin{aligned}
		\frac{\partial u(x,y)}{\partial \bm{n}}  & = \cos(x^2+y), \qquad (x,y) \in \Gamma_1, \\
		u(x,y)  & = \sin(x^2+y), \qquad (x,y) \in \Gamma_2.
	\end{aligned}
\end{eqnarray*}
The analytical solution and the source term are given, respectively, as
\begin{eqnarray*}
	\begin{aligned}
		u(x,y) = \sin(x^2+y), \
		f(x,y) = -2\cos(x^2+y)-(4x^2+1)\sin(x^2+y).
	\end{aligned}
\end{eqnarray*}
We test two { quasi-uniform point distributions with $N=2500$  under two mesh ratio, see Figure~\ref{F:NT2Nodes}}. Nonuniform distributed nodes {with large mesh ratio } are obtained by  the transformation $z = \sin \left( (\pi z)/2 \right)$ in all  coordinates.

\begin{figure}
	\centering
\qquad
	\subcaptionbox{}{
		\includegraphics[width=4.4cm]{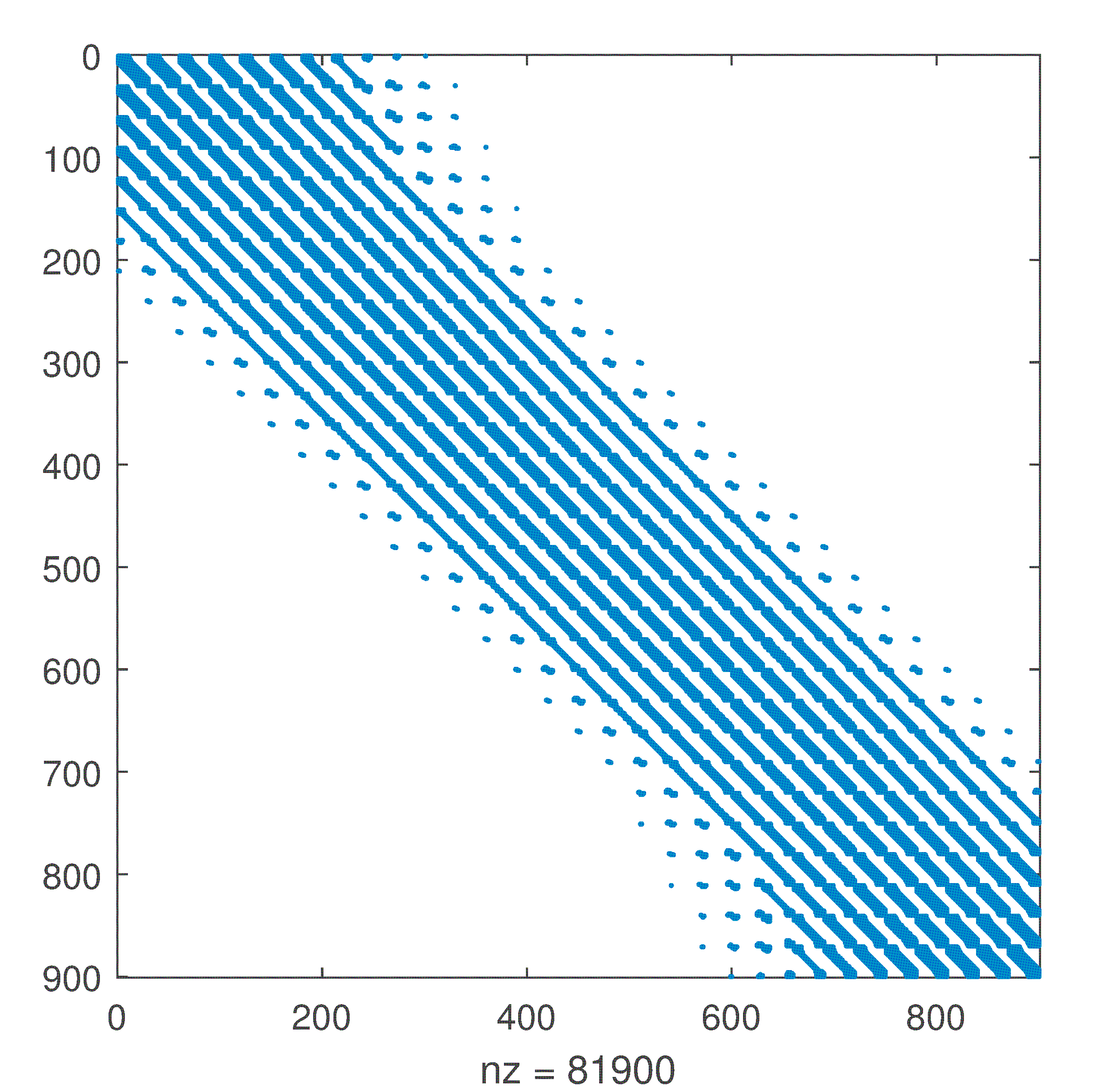}
		\put(-140,30){\rotatebox{90} {\scriptsize Non-adaptive RBF-FD}}
		\put(-130,50){\rotatebox{90} {\scriptsize Row($A$)}}
		\put(-80,-5){ {\scriptsize Col$(A)$}}
	}
	\quad\quad
	\subcaptionbox{}{
		\includegraphics[width=6.2cm]{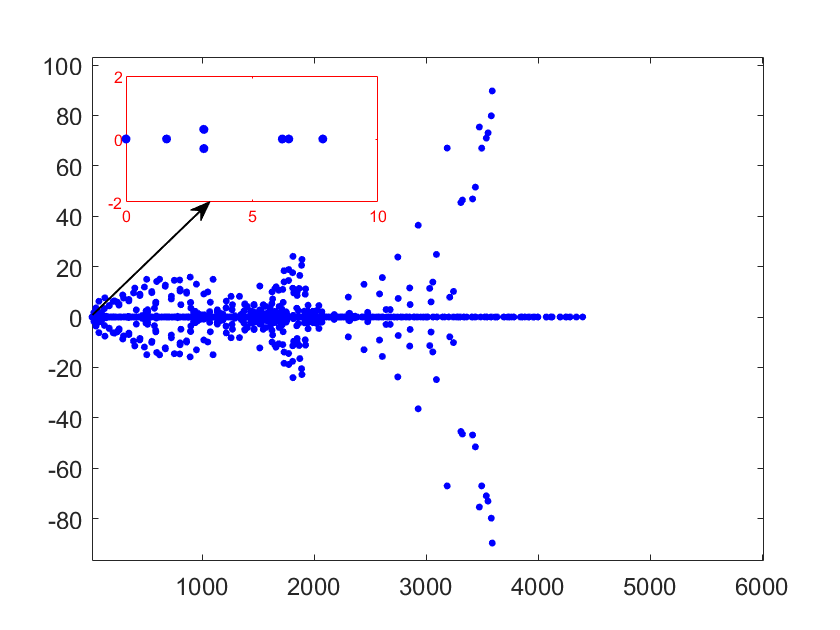}
		\put(-180,50){\rotatebox{90} {\scriptsize Im($\lambda$)}}
		\put(-100,0){ {\scriptsize Re($\lambda$)}}
	}
	\
	\subcaptionbox{}{
		\includegraphics[width=6cm]{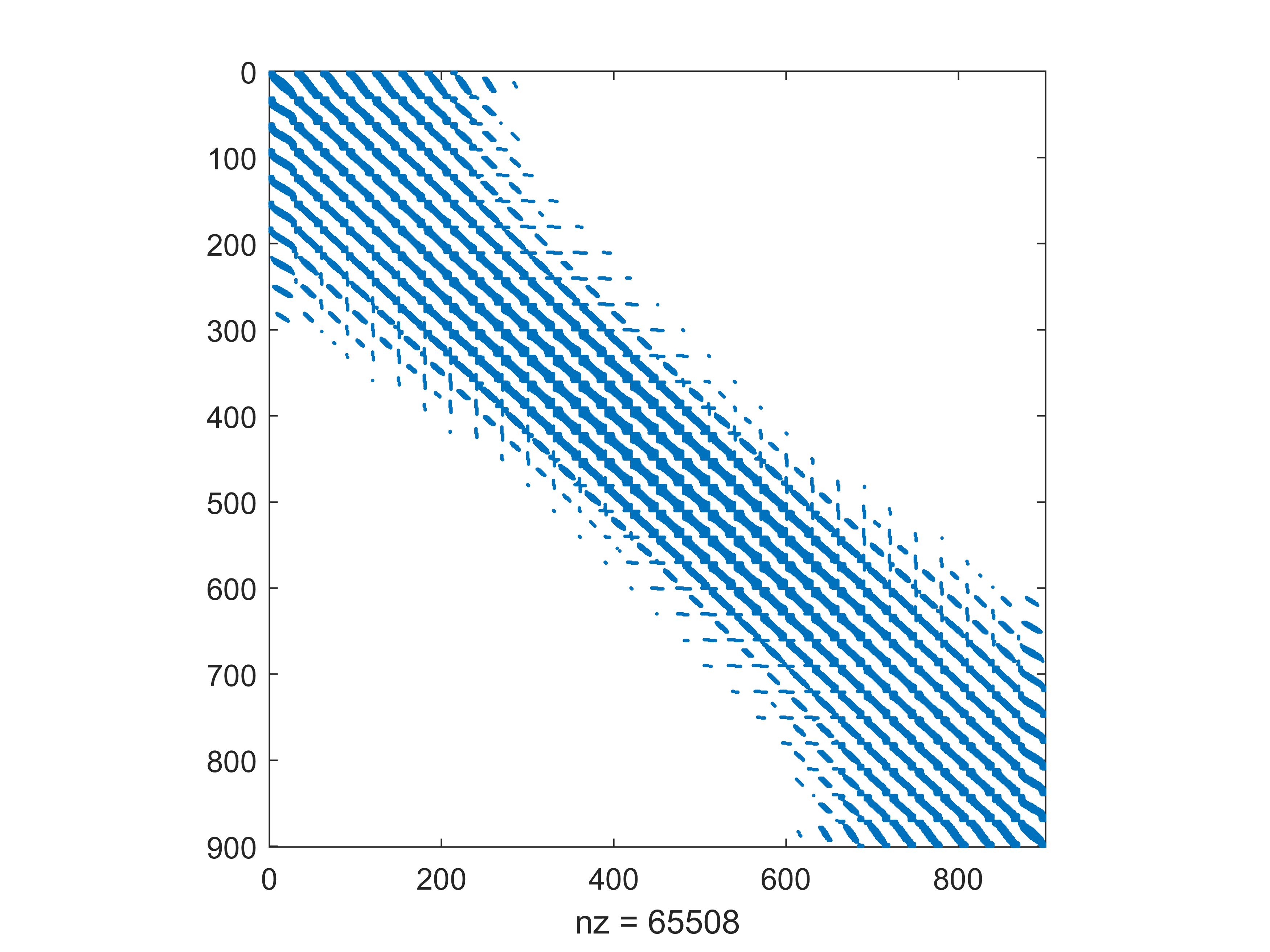}
		\put(-165,30){\rotatebox{90} {\scriptsize Adaptive RBF-FD}}
		\put(-155,50){\rotatebox{90} {\scriptsize Row($A$)}}
		\put(-100,-5){ {\scriptsize Col$(A)$}}
	}
\
	\subcaptionbox{}{
		\includegraphics[width=6cm]{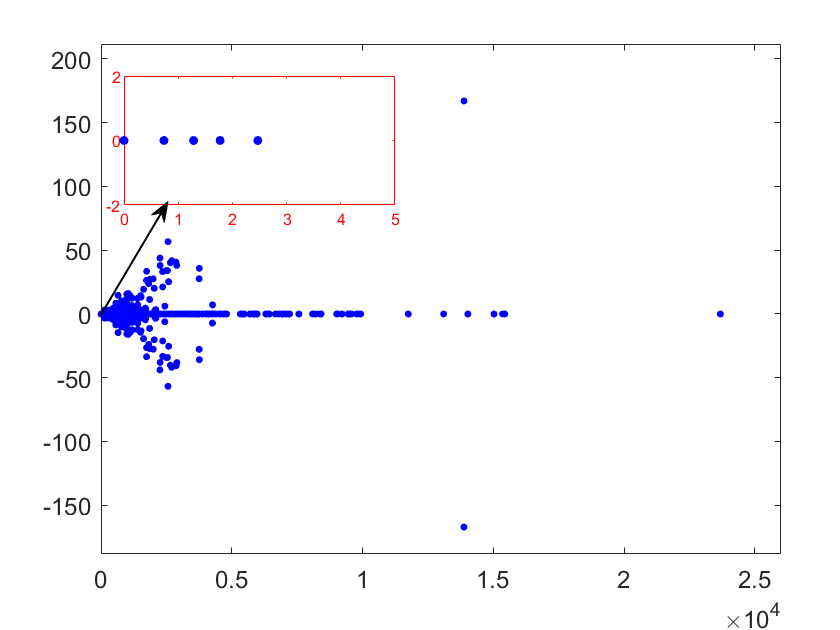}
		\put(-170,50){\rotatebox{90} {\scriptsize Im($\lambda$)}}
		\put(-100,0){ {\scriptsize Re($\lambda$)}}
	}
	\caption{Example 1: results by adaptive PHS+Poly RBF-FD  method for uniform and nonuniform node distribution by $N=900$, (a,c): Sparsity pattern and (b,d): eigenvalue spectra of the system matrix. }
	\label{F:NT2spyeig}
\end{figure}

{For a target convergence rate of $g=3$, Figure~\ref{F:NT2spyeig}~(a-b) show the sparsity patterns and eigenvalues of $\mathcal{L}$ of the globally assembled differentiation matrix for uniform nodes and Figure~\ref{F:NT2spyeig}~(c-d) for nonuniform nodes. Recall that the adaptive method is equivalent to the standard RBF-FD on uniform node distribution. From Figure~\ref{F:NT2spyeig}~(a,c), it can be seen that the adaptive method has a smaller bandwidth  { in the range $[12,15]$ (except near boundaries) compared to  the fixed bandwidth $ [15]$} for standard RBF-FD.  The total nonzero elements are $6.5\times 10^4$, an 20\% reduction.}
	From Figure~\ref{F:NT2spyeig}~(b,d), the real parts of the eigenvalues computed from the differentiation matrix are all non negative, showing that stable approximations can be obtained by our adaptive approach with higher sparsity.

\begin{figure}[t]
	\centering
		\includegraphics[width=7cm]{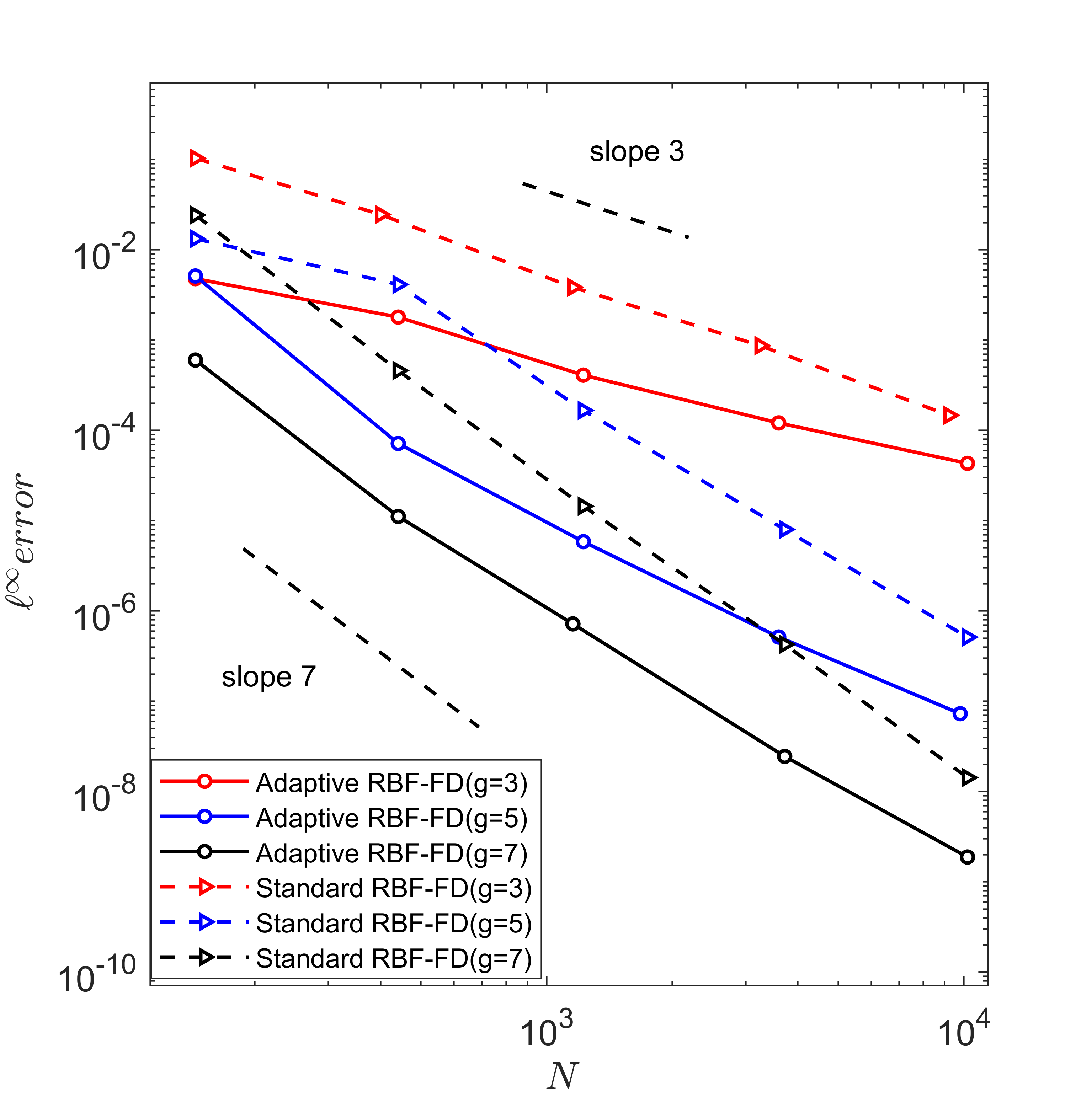}
		\caption{Example 1:  convergence comparison of the adaptive PHS+Poly RBF-FD  algorithm with  nonuniform nodes and standard RBFFD with uniform  nodes }
		\label{F:NT2conv}
	\end{figure}

Figure~\ref{F:NT2conv} shows the convergence profile of $\ell^\infty$-error:
\[
\text{$\ell^{\infty}$-error }=\max _{1 \leq i \leq N}\big|u(x_i)-\hat{u}(x_{i})\big|,
\]
versus $N^{1/d}$
for different target orders of convergence $g=3,5,7$ with $N$ and $d$  being the number of nodes in the domain and dimension of the problem.
The adaptive RBF-FD gives higher accuracy while maintaining a similar (desired) convergence rate and at the same time, it leads to a sparser system matrix.
{For $g=7$ and $N=2500$, the adaptive method has $1.72\times 10^5$ nonzeros versus $2.14\times 10^5$ for standard RBF-FD.}  {From Figure~\ref{F:NT2spyeig} and Figure~\ref{F:NT2conv}, the proposed method is cheaper and more accurate for the example.}


Table~\ref{tab:NT2histo} shows the frequency of polynomial degree when  our adaptive scheme is applied to larger data sets with $N=6241$. For uniform nodes, $p=g+1=8$ were used at all centers. For non-uniform nodes, $p$ ranges from $5$ to $8$ based on local node density.
Hence, the proposed adaptive PHS+Poly RBF-FD  scheme achieves higher accuracy at a lower computational cost by adaptive selection of augmented polynomial degree $p_c$.

	\begin{table}[t]
		\centering
		\small
		\caption{ {Example 1:  the frequency of polynomial degree by the adaptive PHS+Poly RBF-FD  algorithm for uniform and nonuniform node distribution in Fig. \ref{F:NT2Nodes}(b)}}
		\begin{tabular}{ccccc}
			\toprule
			Polynomial degree ($p_i$) & 5 & 6 & 7 & 8 \\
			\midrule
			Non-adaptive RBF-FD   &  0 &  0 & 0 & 6241\\
			\\
			Proposed adaptive RBF-FD & 110 & 1727 & 2694 & 1710 \\
			\bottomrule
		\end{tabular}
		\label{tab:NT2histo}
	\end{table}
	
\begin{figure}
	\centering
	\subcaptionbox{}{
		\includegraphics[width=3.5cm]{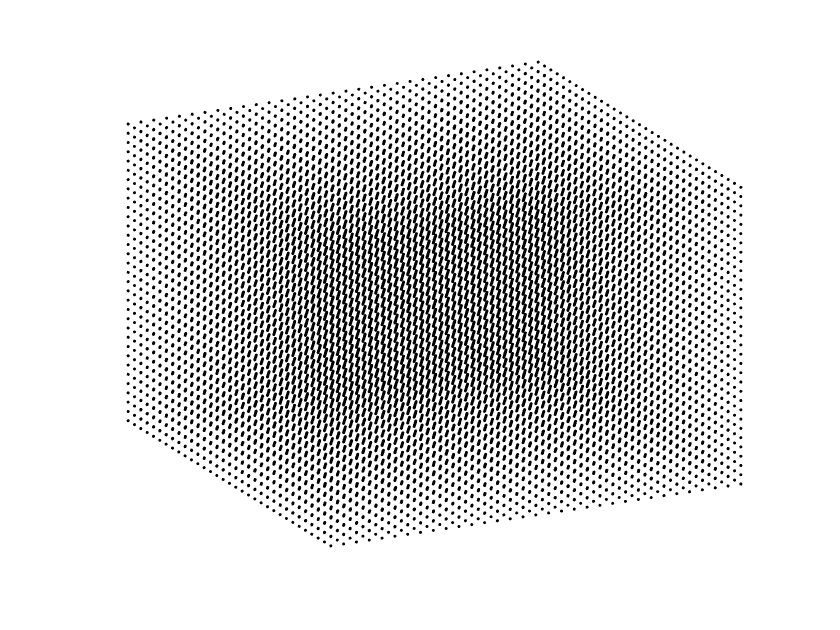}
	}
	\quad
	\subcaptionbox{}{
		\includegraphics[width=3.5cm]{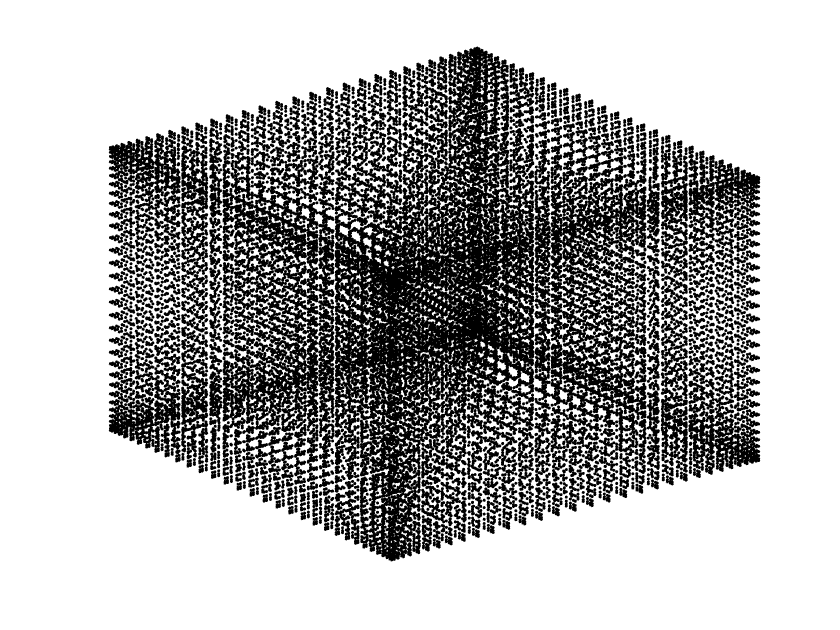}}
	\quad
	\subcaptionbox{}{
		\includegraphics[width=3.5cm]{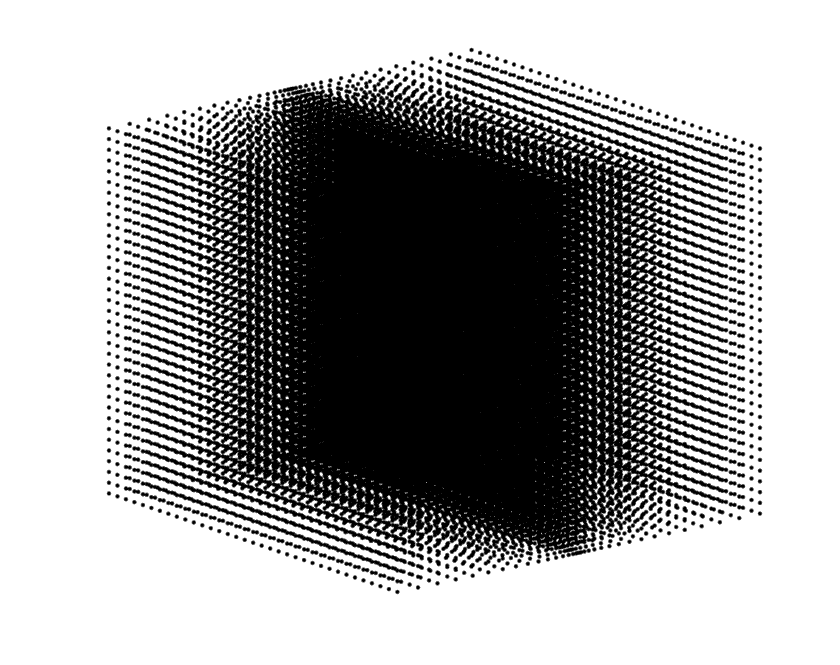}
	}
	\caption{{\normalsize  Example 1: For 3D case, in $\Omega=[-1,1]^3$ with $N=42875$ (a) uniform node distribution,  (b) non-uniform distributed nodes with large mesh ratio generated by sine-transform, and (c) non-uniform nodes with large mesh ratio generated by signed-squared transform at $y$ axis.}}
	\label{F:NT3Nodes3d}
\end{figure}
\begin{figure}[htbp]
	\centering	\includegraphics[width=8cm]{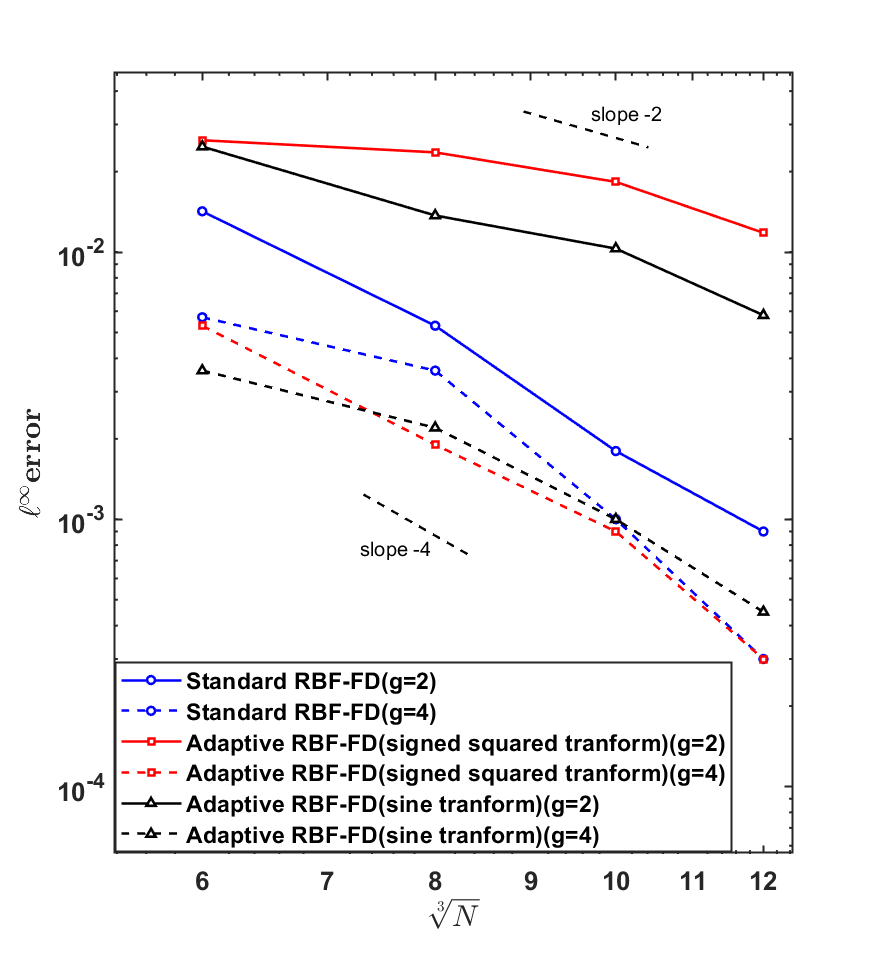}
	\caption{\normalsize {{ Example 1: For  3D case,  convergence comparison of RBF-FD with uniform nodes and the adaptive RBF-FD with non-uniform nodes by $g\in\{2,4\}$.}}}
	\label{F:NT3conv3d}
\end{figure}

{What's more, we also extend  the adaptive method  to solve three dimensional problems. The Poisson equation $-\nabla^2 u(x,y,z)  = f(x,y,z)$ is solved in domain  $\Omega = [-1,1]^3$.  The exact solution is set as $u^{*}(x,y,z)=\sin \left(\pi x / 2\right) \cos \left(\pi y / 2\right) \cos \left(\pi z / 2\right)$. The Dirichlet boundary condition is imposed on the $\partial \Omega$.   Figure~\ref{F:NT3Nodes3d}  shows interior nodes of $N=42875$ by uniform points (a) and non-uniform points (b) and (c) which generated by sine-transformation in all directions and sign-transformation  $z=sign(z)*z^2$ along $y$ axis respectively.
 Figure~\ref{F:NT3conv3d} shows the convergence behavior of the RBF-FD with uniform nodes and the adaptive RBF-FD with non-uniform nodes shown in Figure~\ref{F:NT3Nodes3d} by setting the global order $g\in \{2, 4\}$.  It can be seen that both the adaptive RBF-FD and standard RBF-FD can obtain desired convergence order in 3D cases.  For the global convergence rate $g=4$, Table~\ref{tab:NT3histo} presents frequency of the three-dimensional polynomial degree $(p_i)$ by the standard  and the adaptive RBF-FD methods under two kinds of large mesh ratio node distribution  respectively.  It shows again that the adaptive method can have stable solutions with low computational cost  by adaptively selecting the augmented polynomial basis. }

{
\begin{table}
	\centering
	\small
	\caption{\normalsize {Example 1: For  3D case,  fiixed convergence rate $g=4$ and node number $N=13824$,  the frequency of  polynomial degree $(p_i)$ by the adaptive RBF-FD algorithm for uniform and non-uniform node distributions in Figure 7.}}
 \begin{tabular}{ccccc}
	\toprule
		 {Polynomial degree} ($p_i$) & 3 & 4 & 5  \\
		\midrule
		RBF-FD   &  0 &  0 & 13824 \\
		\\
		Adaptive RBF-FD(sine transform) & 800 & 11984 & 1040 \\
		\\
		Adaptive RBF-FD(signed squared transform) & 0 &5776 &8048\\
		\bottomrule
	\end{tabular}
	\label{tab:NT3histo}
\end{table}	}

	\subsection{Example 2: Benchmark test for  adaptive FEM}
\begin{figure}[t]
		\centering
		\includegraphics[width=8cm]{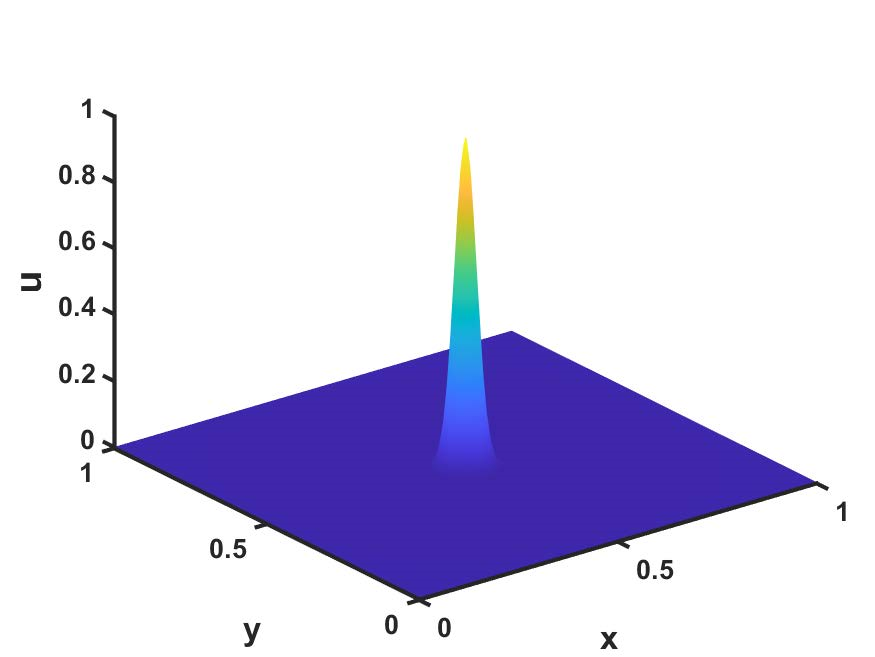}
		\caption{ Example 2: A reference solution with $\alpha = 1000$ and $(x_c,y_c) = (0.5,0.5)$}
		\label{F:NT3exact}
	\end{figure}

	This example considers a benchmark test problem  from  US National Institute for Standards and Technology (NIST)  \cite{mitchell2013collection} for   adaptive finite element method algorithms.
	It is a Poisson equation
	$
	-\Delta u(x,y) - f(x,y) = 0
	$ in $ \Omega = [0,1]^2$ with
	exact solution being an exponential peak:
	\begin{equation*}
		u(x,y) =  \exp\big(-\alpha  ( (x-x_0)^2+(y-y_0)^2 )\big),
	\end{equation*}
	where $(x_0,y_0)$ is the peak location  and $\alpha$ controls  the peak strength.  Dirichlet boundary conditions matching the exact solution are imposed. A typical value $\alpha = 1000$ is suggested for the test. The right hand side function $f(x,y)$ satisfying the exact solution is
	\begin{equation*}
		f(x,y) = {-}4 \left( \alpha^2  (x-x_0)^2+ \alpha^2  (y-y_0)^2 \right)\exp\big(-\alpha  ( (x-x_0)^2+(y-y_0)^2 )\big) .
	\end{equation*}
	The reference solution with $\alpha = 1000$ is shown in the Figure~\ref{F:NT3exact}.

\begin{figure}
		\centering
		\subcaptionbox{}{
			\includegraphics[width=4.4cm]{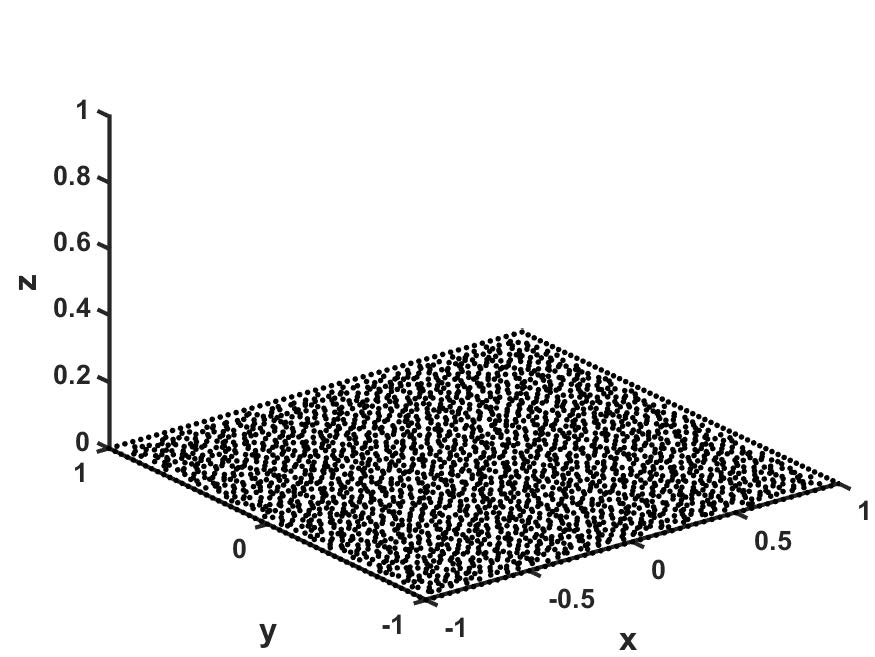}
			\put(-130,20){\rotatebox{90} {\scriptsize Standard RBF-FD}}
		}
		\subcaptionbox{}{
			\includegraphics[width=4.4cm]{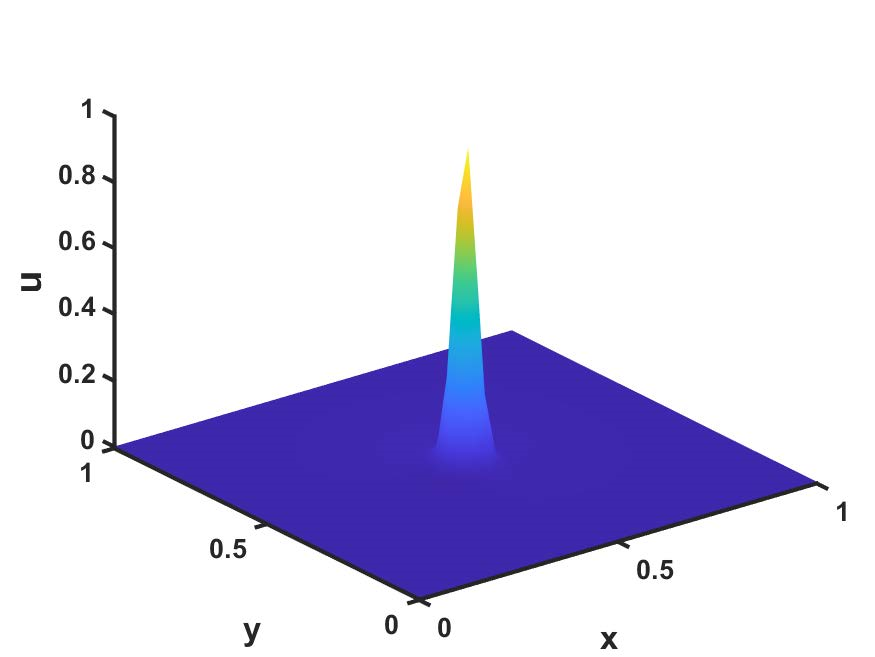}
		}
		\subcaptionbox{}{
			\includegraphics[width=4.4cm]{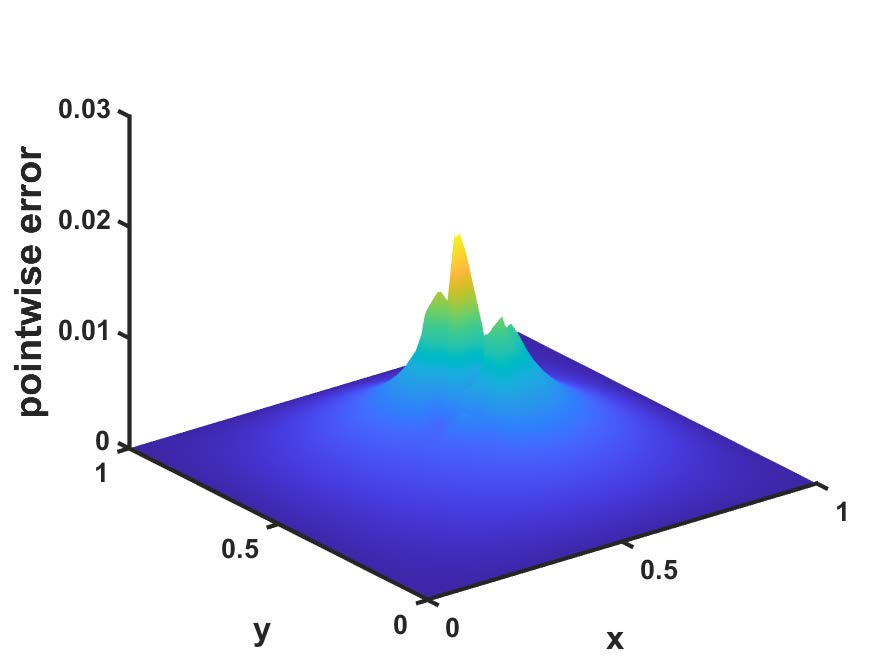}
		}
		\\
		\subcaptionbox{}{
			\put(-8,20){\rotatebox{90} {\scriptsize Adaptive RBF-FD}}
			\includegraphics[width=4.4cm]{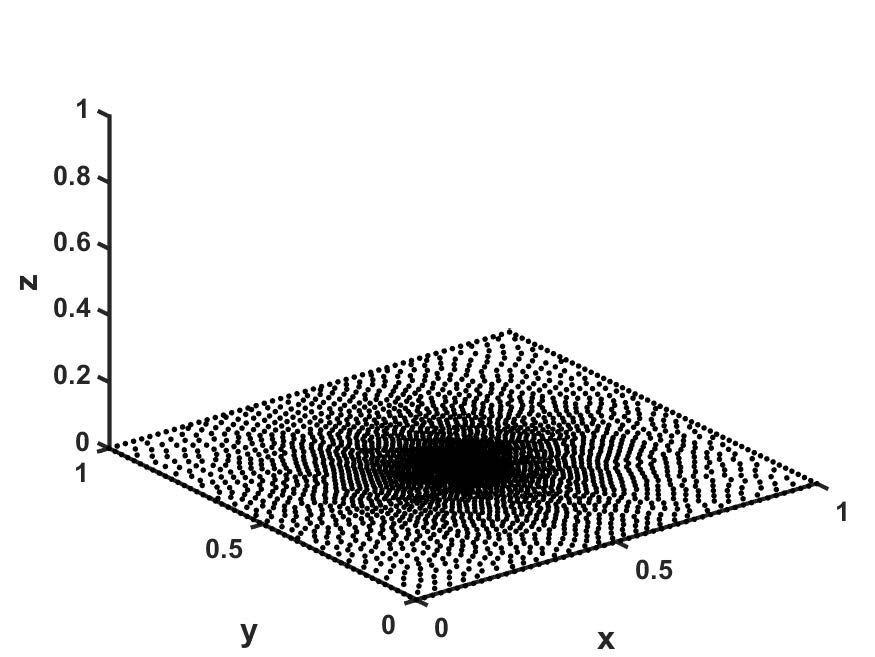}}
		\subcaptionbox{}{
			\includegraphics[width=4.4cm]{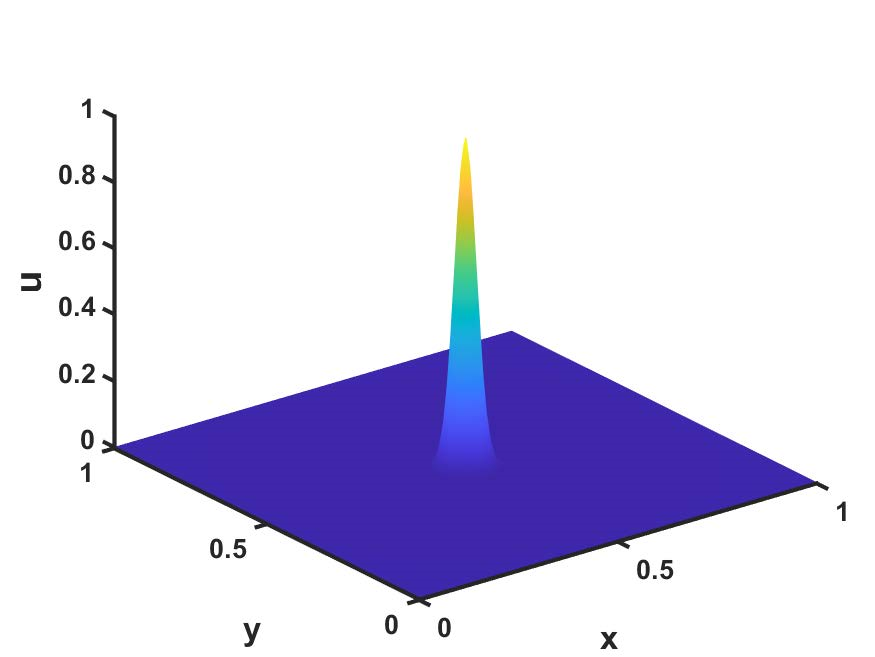}}
		\subcaptionbox{}{
			\includegraphics[width=4.4cm]{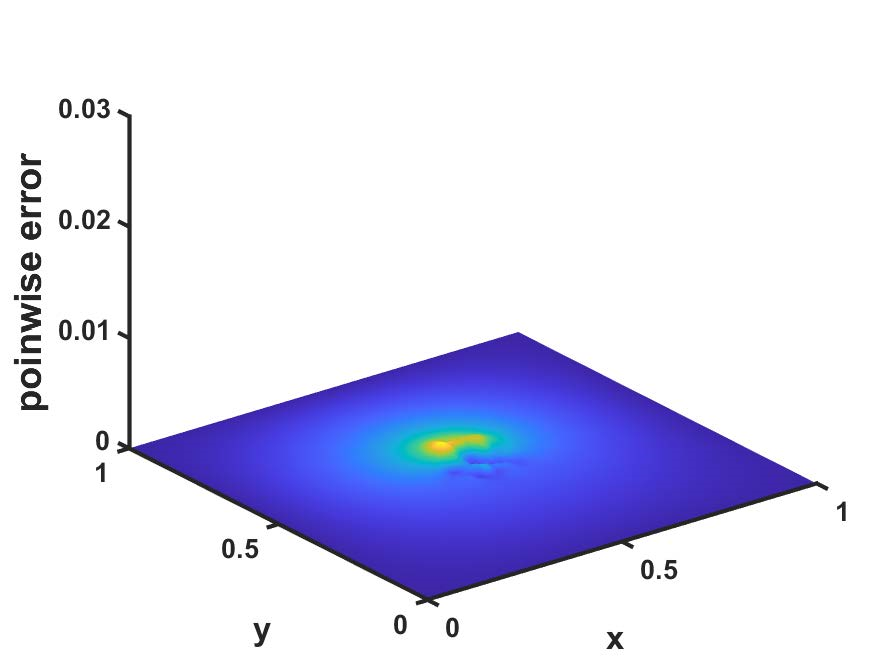}}
		\caption{ {Example 2: setting $N = 2470, \ g=5$, results by adaptive RBF-FD approach under uniform node-distribution with  $p=6$ in first row and  under nonuniform nodes with variable polynomial degree in second row, (a, d) node distribution, (b, e) approximated solution, and (c, f) pointwise  absolute error} }
		\label{F:NT3sol}
	\end{figure}

	Firstly, we solve the problem with standard RBF-FD approach over quasi-uniform nodes in the domain and a fixed polynomial degree $p=6$. {Figures~\ref{F:NT3sol}~(a,b,c)}  show results for a quasi-uniform set of $N=2470$ nodes: nodes, numerical solution, and  error function.
	The error is largest near the peak, indicating inaccuracy there.
	
\begin{figure}
		\centering
			\includegraphics[width=7cm]{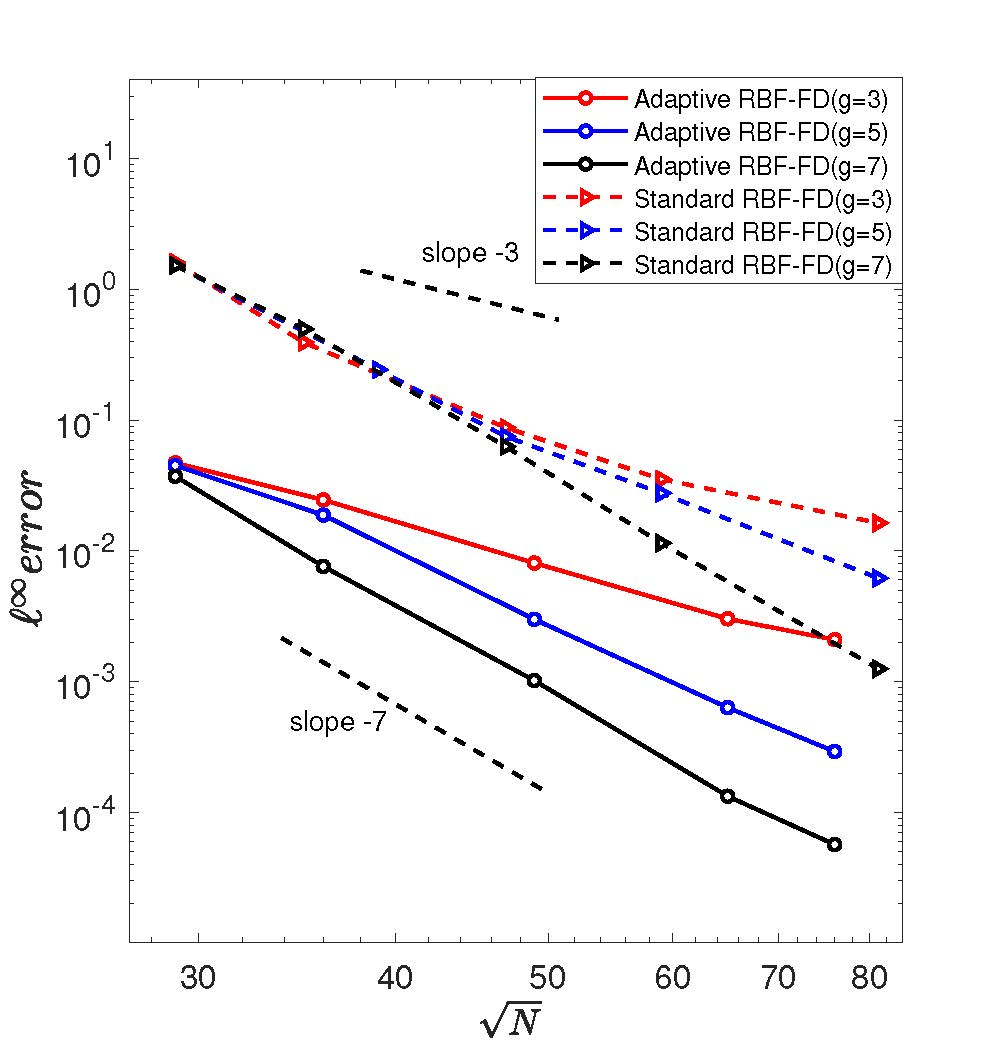}
		\caption{ Example 2: Comparison of convergence of the standard RBF-FD and the Adaptive RBF-FD for a Poisson equation with delta-like exact solution.}
		\label{F:NT3conv}
	\end{figure}

	We use { the \texttt{NodeLab} algorithm \texttt{NodeLab}  \cite{Mishra2019nodelab,FORNBERG2015531} }to generate data-dependent nodes with same number $N=2470$, which is denser in the peak zone, see
	Figures~\ref{F:NT3sol}~(d).  This adaptive node-distribution enhance the accuracy in the approximation. {Figures~\ref{F:NT3sol}~(e,f)} show the approximated solution and point-wise absolute error with an expected global order $g=5$. Improved in accuracy around the peak zone is obvious in comparison  with standard RBF-FD. By setting convergence rate $g\in \{3,5,7\}$, Figure \ref{F:NT3conv} shows the $g^{th}$ convergence rate can be obtained by both method, but the proposed adaptive method are more accurate by more than one order of magnitude.
	
	\subsection{Example 3: Solving a heat equation}
	We consider a time-dependent heat equation incorporated with compatible   initial and boundary conditions:
	\begin{equation*}
		\left\{
		\begin{aligned}
			u_t- \bigtriangleup u &= f(x,y,t) \qquad  \ \text{for} \ (x,y,t) \in \Omega \times (0,T],\\
			u(x,y,0)&=g_1(x,y)  \quad \qquad \text{for} \ (x,y) \in \Omega,\\
			u(x,y,t)&=g_2(x,y,t) \qquad \ \text{for} \ (x,y,t) \in \partial \Omega \times (0,T].
		\end{aligned}
		\right.
		\label{eq:heat}
	\end{equation*}
	in $\Omega =[-1,1]^2$ with exact solution  $u^*(x,y,t)=\exp(t)\cos(\pi x)\cos(\pi y)$.
	The initial condition and boundary conditions are generated from the exact solution:
	\[
	g_1(x,y)=\cos(\pi x)\cos(\pi y), \qquad g_2(x,y,t) = \exp(t)\cos(\pi x)\cos(\pi y),
	\]
	and so as the source term  $f(x,y,t)=(1+2\pi^2)u^*(x,y,t)$.

	We discretize the PDE \eqref{eq:heat} in time using the $\theta$-method over $M$ time steps with step size $\tau=T/M$:
	\begin{equation*}
		\begin{split}
			\frac{u^n-u^{n-1}}{\tau}=\theta (\bigtriangleup u^n+f^n)+(1-\theta)(\bigtriangleup u^{n-1}+f^{n-1}).
		\end{split}
	\end{equation*}
	We consider both the  backward-Euler scheme with $\theta=1$:
	\begin{equation*}
		\begin{split}
			u^n- \tau \bigtriangleup u^n &=u^{n-1}+\tau f^n(x,y,t),
		\end{split}
		\label{eq:timeEuler}
	\end{equation*}
	and   Crank-Nicolson scheme with   $\theta=1/2$:
	\begin{equation*}
		\begin{split}
			u^n- \frac{\tau}{2} \bigtriangleup u^n &=u^{n-1}+\frac{\tau}{2} \bigtriangleup u^{n-1}+\frac{\tau}{2} \bigtriangleup (f^n+f^{n-1}).
		\end{split}
		\label{eq:timeCN}
	\end{equation*}
	
\begin{figure}[htbp]
		\centering
		\setlength{\abovecaptionskip}{0pt}
		\setlength{\belowcaptionskip}{0pt}
		\begin{tabular}{ccccc}
			\begin{overpic}[height=5.5cm]{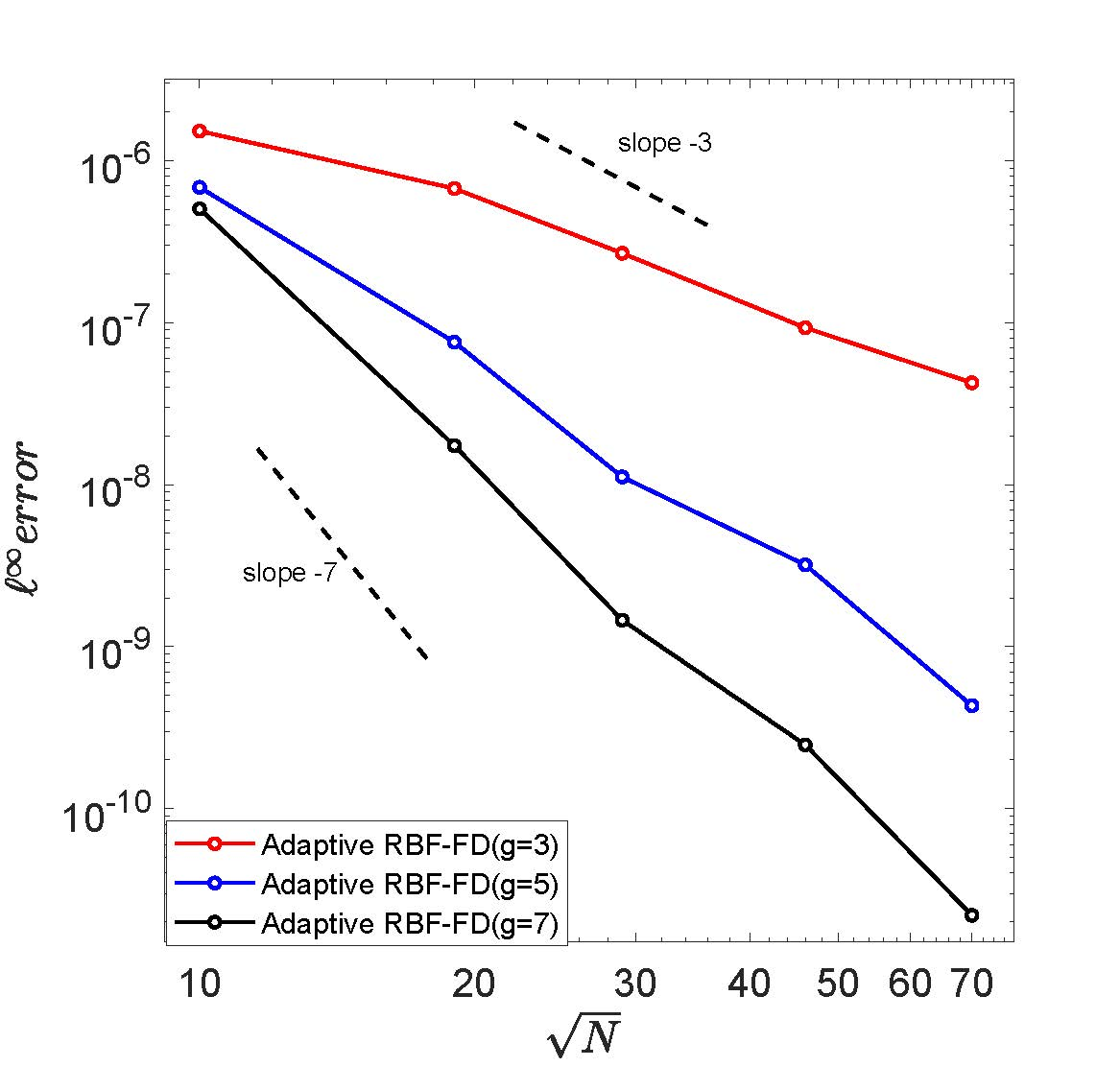}
				\put(47,90) {\small $(a)$}
			\end{overpic}
			&
			\begin{overpic}[height=5.5cm]{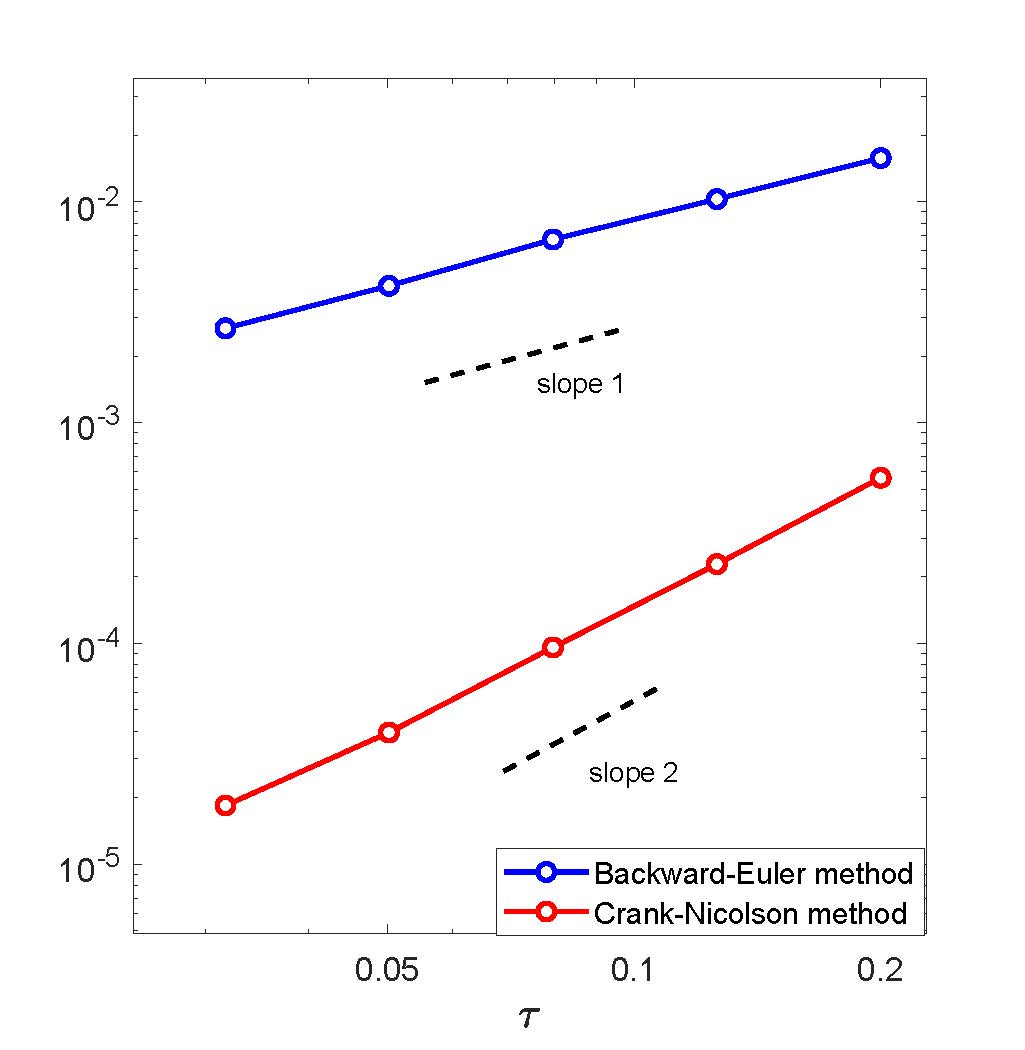}
				\put(45,94) {\small $(b)$}
			\end{overpic}
		\end{tabular}
		\caption{{Example 3: the convergence results of adaptive RBFFD (PHS+Poly) for solving the heat equation (\ref{eq:heat}): (a)
		space convergence by $\tau=1E-7, T=1E-6,\ N\in [10^2, 70^2]$; (b) time convergence by $\tau\in [0.02,0.2]$, $T=1$, $N=64^2$}}
		\label{fig:ex3conv}
	\end{figure}

	The nonuniform discrete points in the domain is generated by the same method as described in Example 1.  The user-defined global convergence rate set to $g=\{3,5,7\}$ and  {PHS kernel is $\phi(r)=r^{2m-1}, m=3$}. Figure~\ref{fig:ex3conv} illustrate the convergence behavior of the solutions by the $\ell^\infty(\Omega,T)$ error. Figure~\ref{fig:ex3conv}~(a) plots the space convergence results with respect to $\sqrt{N}$ by $\tau=1E-7, T=1E-6$, discrete points in domain nodes $N\in [10^2,70^2]$.
	 The local stencil is set as $n_s=100$. The result shows that the convergence rate $g$ is obtained. Figure~\ref{fig:ex3conv}~(b) illustrates the time convergence results by $T=1$, $N=64^2$. It can be observed that the Crank-Nicolson method achieves second-order accuracy in time, while the backward-Euler method achieves first-order accuracy in time. These results suggest that the proposed adaptivity does not cause instability when solving parabolic equations.
	
	\subsection{Example 4: Application in elastic wave model}
	{We apply the adaptive RBF-FD method to the following elastic wave Navier equation \cite{yuan2019trefftz}:
		\begin{equation}
			\left\{\begin{array}{l}
				\mu \Delta u_1 + (\lambda+\mu)({\partial_x^2 u_1}+{\partial_x\partial_y u_2})+\omega ^2 \rho u_1 =f_1\\
				\mu \Delta u_2 + (\lambda+\mu)({\partial_x\partial_y u_1}+ {{\partial_y^2} u_2} )+\omega ^2 \rho u_2=f_2 \\
			\end{array},\right.
			\label{Navier}
		\end{equation}
with $\omega$ being the angular frequency and $\lambda,\ \mu$ being Lam{\'e} constants. We assume the medium is homogeneous  and set the density $\rho \equiv 1$ after normalization. The time-harmonic displacement vector $\boldsymbol{u}=(u_1,u_2)^T$ satisfies the lowest-order absorbing boundary condition:
		\begin{equation*}
			\mathbf{T}^{(\mathbf{n})}(\boldsymbol{u})-i \eta \boldsymbol{u}=\mathbf{g} \quad \text { on } \ \gamma=\partial \Omega,
			\label{Nav_bc}
		\end{equation*}
		where the traction operator $\mathbf{T}^{(\boldsymbol{n})}$ is
		\begin{equation*}
			\mathbf{T}^{(\boldsymbol{n})}(\boldsymbol{u})=2 \mu \frac{\partial \boldsymbol{u}}{\partial \boldsymbol{n}}+\lambda \boldsymbol{n} \nabla \cdot \boldsymbol{u}+\mu \boldsymbol{n} \times(\nabla \times \boldsymbol{u}),
			\label{Tn}
		\end{equation*}
		with $\boldsymbol{n}$ being the unit normal vector and a real-valued matrix function $\eta$ is defined by
		{\begin{equation*}
			\begin{aligned}
				\eta=\omega \rho\left(C_{\mathrm{P}} \boldsymbol{n} \boldsymbol{n}^{T}+C_{\mathrm{S}} \boldsymbol{s}  \boldsymbol{s}^{T}\right) ,
				C_{\mathrm{P}}=\sqrt{\frac{\lambda+2 \mu}{\rho}} , \ C_{\mathrm{S}}=\sqrt{\frac{\mu}{\rho}}
			\end{aligned}
			\label{eta}
		\end{equation*}}
		with $\boldsymbol{s}$ being a unit tangent vector. 
We study the problem in  $\Omega=[0,1]^2$, and nonuniform distribution nodes are obtained by using the transformation $z=(\sin(\pi z)+1)/2$ on each coordinate. The components of the exact displacement field $\mathbf{u}=(u_1,u_2)^T$ is given as follows
		\begin{equation}
			\left\{\begin{array}{l}
				\displaystyle
				u_{1}=\alpha_{\mathrm{S}}\left\{\exp(-\alpha_{\mathrm{S}}y)-\frac{2 \kappa_{R}^{2}}{\kappa_{R}^{2}+\alpha_{\mathrm{S}}^{2}} \exp(-\alpha_{\mathrm{P}} y)\right\} \exp(i \kappa_{R} x) \\
				\displaystyle
				u_{2}=i \kappa_{R}\left\{ \exp(-\alpha_{\mathrm{S}} y)-\frac{2 \alpha_{\mathrm{P}} \alpha_{\mathrm{S}}}{\kappa_{R}^{2}+\alpha_{\mathrm{S}}^{2}} \exp(-\alpha_{\mathrm{P}} y)\right\} \exp(i \kappa_{R} x)
			\end{array}, \right.
			\label{q4_exa}
		\end{equation}
		where $\alpha_{\mathrm{P}}=\sqrt{\kappa_{R}^{2}-\kappa_{\mathrm{P}}^{2}}, \ \alpha_{\mathrm{S}}=\sqrt{\kappa_{R}^{2}-\kappa_{S}^{2}}$.
		We take the Lam{\'e} constants $\lambda=2$, $\mu=1$, angular frequency $\omega=2\pi$, which decide the value of $\kappa_R=2.14\pi, \ \kappa_P=\pi, \ \kappa_S=2\pi$. We report the relative $L^2$-error:
		\begin{equation*}
			E_{\eta}=\frac{\left\|\boldsymbol{u}^{e}-\boldsymbol{u}^{n}\right\|_{2}}{\left\|\boldsymbol{u}^{e}\right\|_{2}}
			=\left({\frac{\int_{\Omega}\left(\overline{\boldsymbol{u}}^{e}-\overline{\boldsymbol{u}}^{n}\right)^{T}\left(\boldsymbol{u}^{e}-\boldsymbol{u}^{n}\right) d \Omega}{\int_{\Omega}\left(\overline{\boldsymbol{u}}^{e}\right)^{T}\left(\boldsymbol{u}^{e}\right) d \Omega}}\right)^{1/2}.
		\end{equation*}	
		Setting the global convergence order  $g\in \{3,5,7\}$,  Figure~\ref{fig:q4}~(a) compares the proposed adaptive PHS+Poly RBF-FD on nonuniform nodes and standard RBF-FD on uniform nodes. The adaptive method achieves higher accuracy while maintaining the desired convergence rate.
		
		Angular frequency $\omega$ determines oscillatory behavior: higher $\omega$ causes more frequent oscillations. Figure~\ref{fig:q4}~(b) plots error $E_\eta$ versus $\omega$ for $g=5$. Errors progressively worsen with increasing $\omega$, but the adaptive RBF-FD is less sensitive.
		Figure~\ref{fig:q4_2} plots the real and imaginary parts of $u_1$,$u_2$ at circle $r=0.5$ for $N=3600$. Both methods recover the analytic solutions well.
		For $g=5$, standard RBF-FD has $8.2\times 10^5$ nonzeros while adaptive RBF-FD has $6.5\times 10^5$, a reduction of $1.7\times 10^5$ nonzeros. For $g=7$, adaptive RBF-FD reduces nonzeros by $2.3\times 10^5$. The adaptive method saves more computation cost at higher accuracy.
For realistic elastic wave problems, the adaptive RBF-FD could be an efficient alternative for practitioners.

		\begin{figure}
			\centering
			\subcaptionbox{}{
				\includegraphics[height=5.2cm]{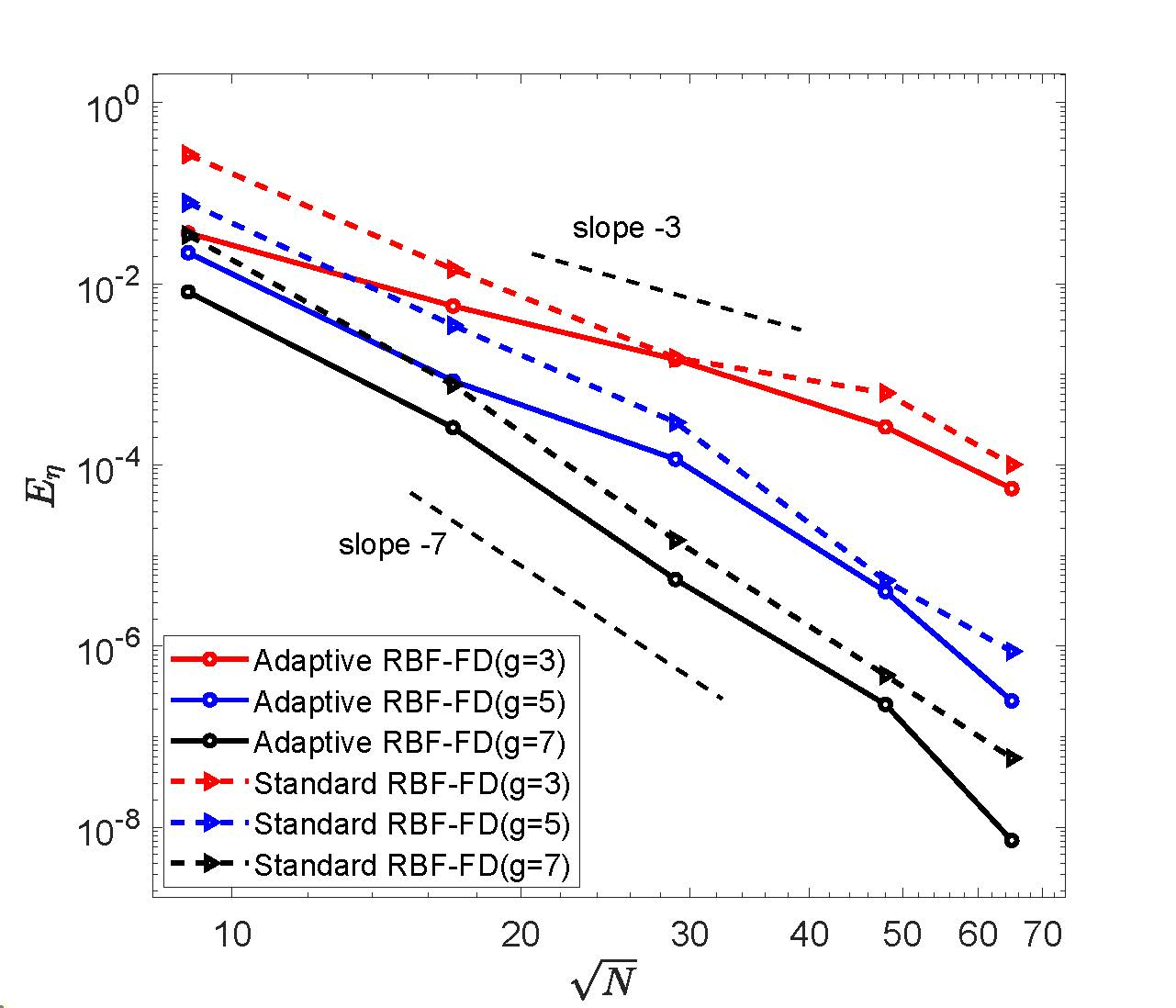}
			}
			\subcaptionbox{}{
				\includegraphics[height=5.2cm]{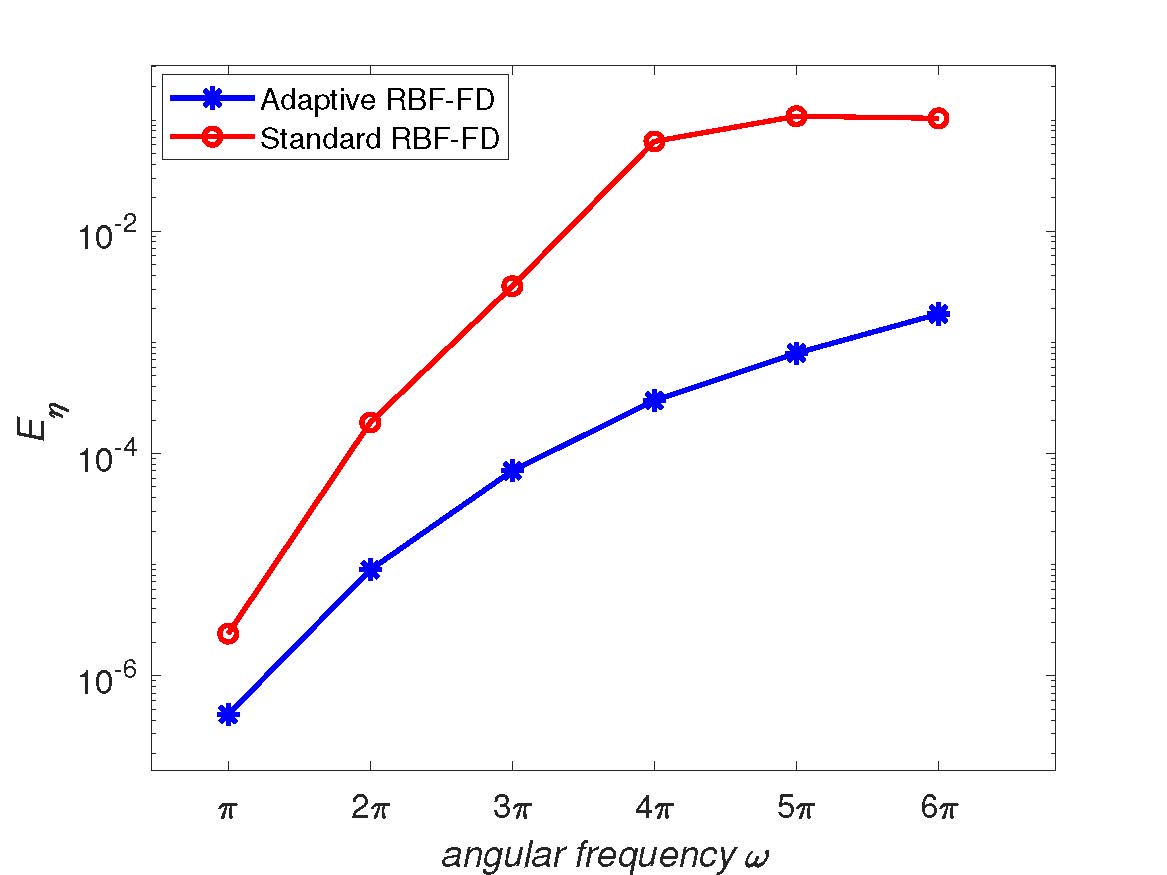}
			}
			\caption{ Example 4: (a) Comparison of convergence of the standard RBF-FD and the Adaptive RBF-FD for Navier equation. (b) The relative errors with respect the angular frequencies.}
			\label{fig:q4}
		\end{figure}
		
		\begin{figure}
			\centering
			\subcaptionbox{}{
				\includegraphics[width=5.5cm]{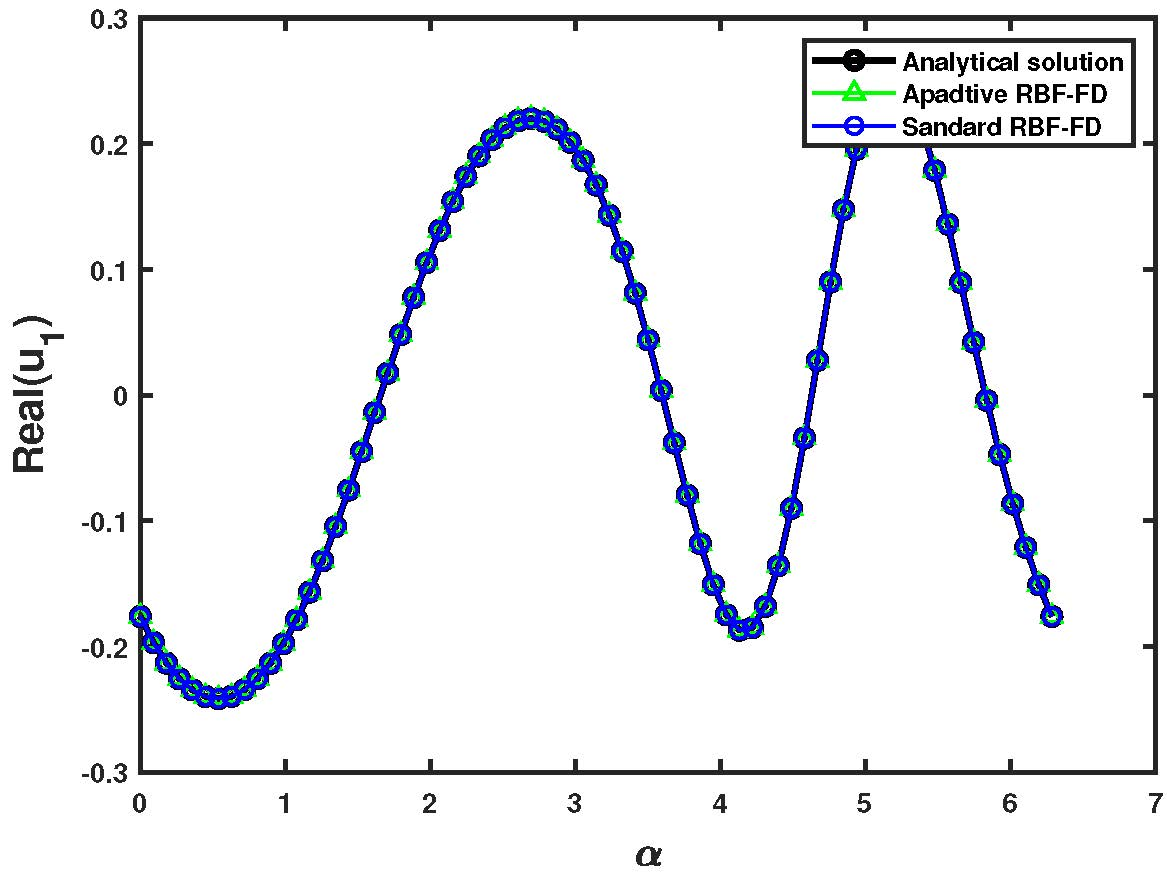}
			}
			\subcaptionbox{}{
				\includegraphics[width=5.4cm]{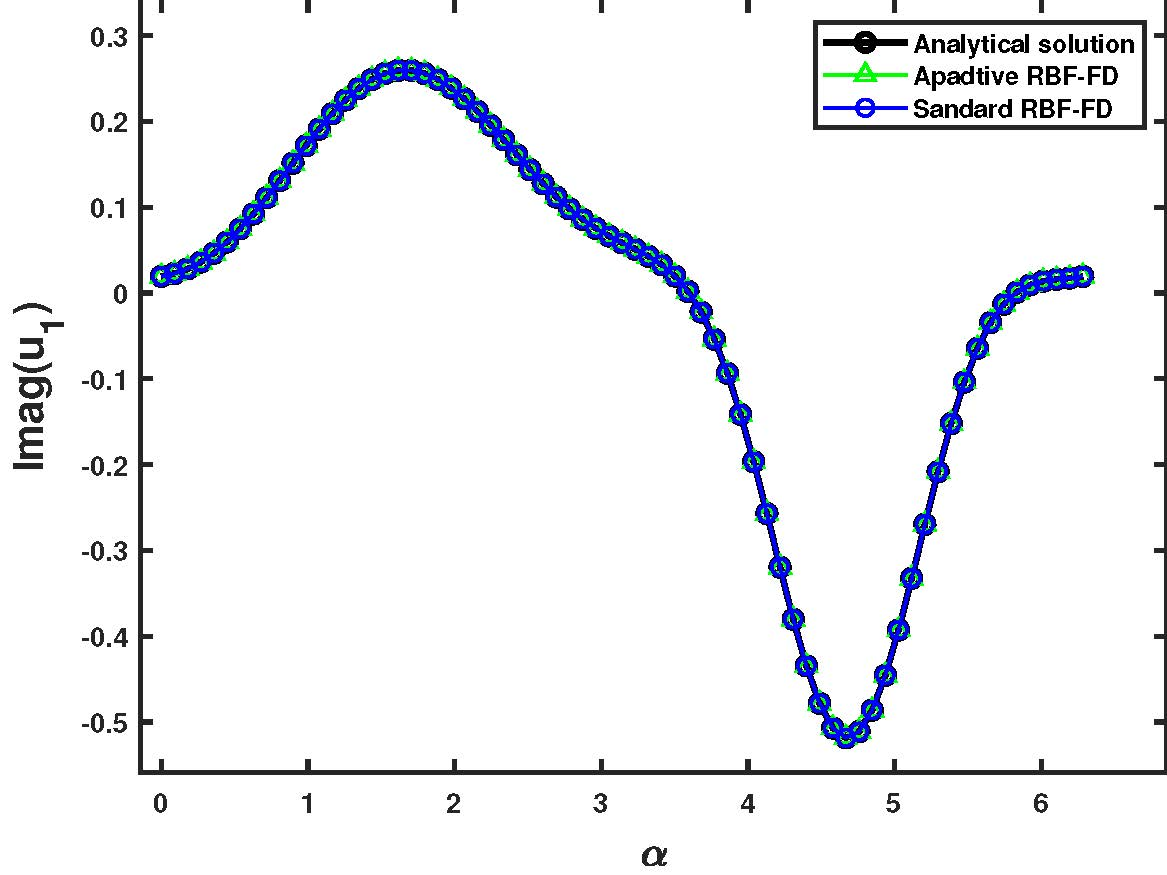}
			}
			\subcaptionbox{}{
				\includegraphics[width=5.5cm]{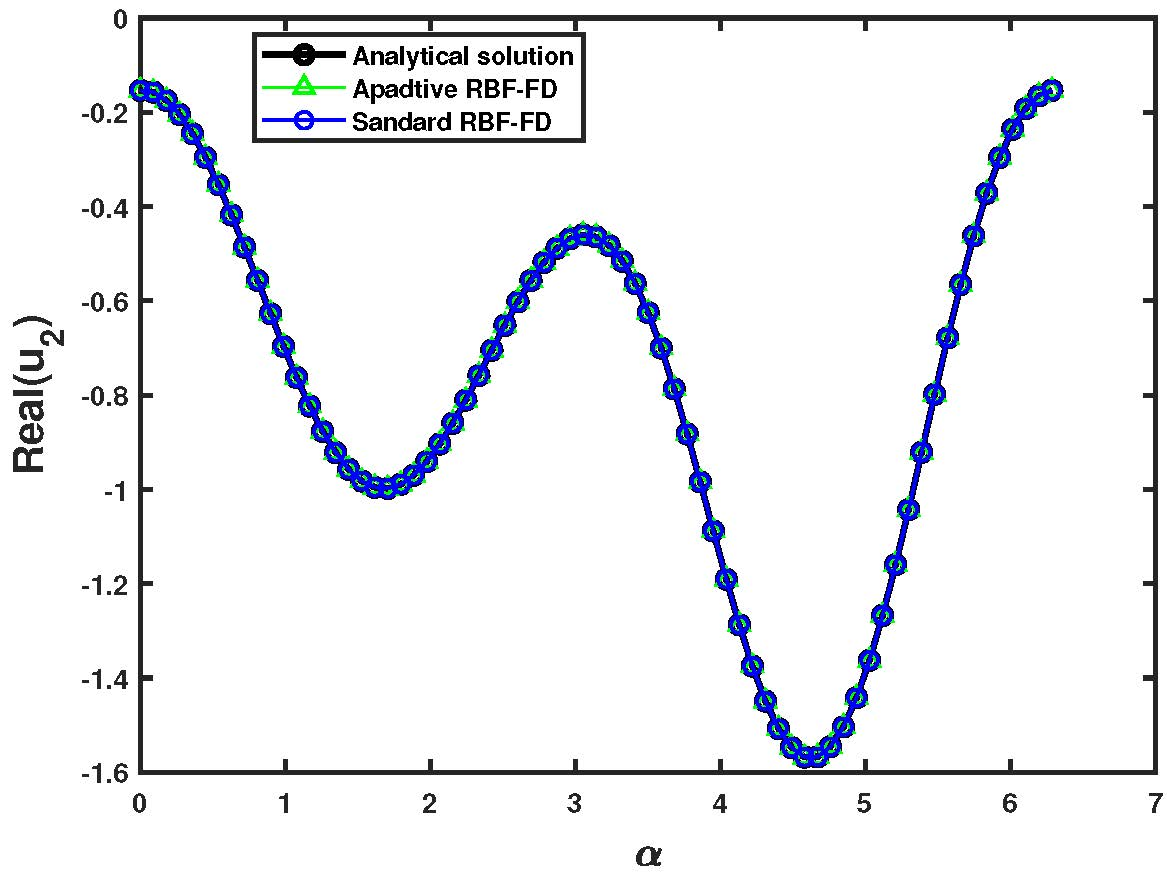}
			}
			\subcaptionbox{}{
				\includegraphics[width=5.5cm]{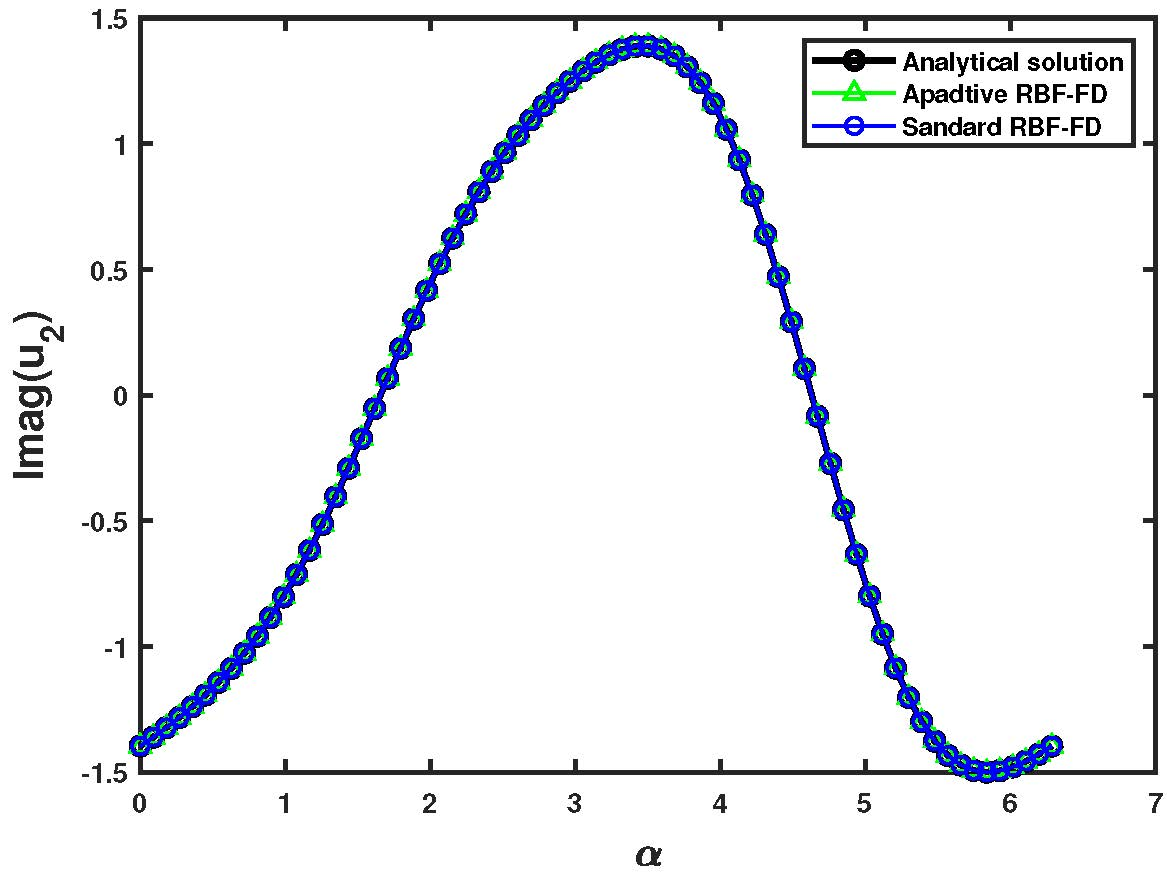}
			}
			\caption{ Example 4:  results  by adaptive RBF-FD for the Navier
				equation with $\omega=2\pi$ and $N=3600$.}
			\label{fig:q4_2}
		\end{figure}
	
		\section{Conclusions}\label{sec:con}
		We developed an adaptive PHS+Poly RBF-FD method well suited for high-order approximations over nonuniform nodes. For uniform nodes, adaptive RBF-FD simplifies to standard RBF-FD. By varying the polynomial degree at each point based on the desired global convergence order, we constructed the adaptive scheme. Our adaptive method provides several advantages: it achieves the specified global convergence order and high accuracy for nonuniform nodes while performing comparably to standard RBF-FD for uniform nodes. The differentiation matrix is sparser, suggesting lower computational cost that will increase for large-scale problems solved using sparse linear solvers. 
{By adaptively choosing local stencil sizes and polynomial degrees based on node density, the proposed method saves the computational resources by  balancing accuracy and efficiency.   In other words,  the proposed method is not aimed to improve accuracy. }  The method also applies to both steady and time-dependent problems, as demonstrated for the heat equation and elastic wave problems. These characteristics make the adaptive RBF-FD method a powerful tool for solving PDEs using scattered data. Future work will explore applications to large-scale and 3D problems to fully leverage the computational efficiency gains of the adaptive method.

		\section*{Acknowledgments}
		This work was funded by the Hong Kong Research Grant Council GRF Grants (12301520,\\ 12301021,12300922), a National Science Foundation of China  (12201449).
		
		\newpage
		\normalem
		\bibliographystyle{siam}
		\bibliography{sample}
		
	\end{document}